\documentclass[12pt,a4paper]{article}
\usepackage[]{inputenc}
\usepackage{amsmath}
\usepackage{float}
\usepackage{multirow}
\usepackage{array}
\usepackage{amsfonts}
\usepackage{amssymb}
\usepackage{color}
\usepackage{authblk}
\usepackage{graphicx}
\usepackage[labelfont=bf]{caption}
\usepackage{caption,subcaption}
\usepackage{hyperref}
\usepackage[left=2cm,right=2cm,top=2cm,bottom=2cm]{geometry}

\begin{document}

\title{Reduced basis stabilization for the unsteady Stokes and Navier-Stokes equations}

\author[1, 2]{\small Shafqat Ali\thanks{ali.qau1987@gmail.com}}
\author[1,3]{\small Francesco Ballarin\thanks{francesco.ballarin@sissa.it}}
\author[1]{\small Gianluigi Rozza\thanks{gianluigi.rozza@sissa.it (corresponding author)}}
\affil[1]{\footnotesize mathLab, Mathematics area, SISSA, via Bonomea 265, I-34136, Trieste, Italy}
\affil[2]{\footnotesize Faculty of Engineering Sciences, Ghulam Ishaq Khan Institute of Engineering Sciences and Technology, Topi, 44000 Pakistan}
\affil[3]{\footnotesize Department of Mathematics and Physics, Catholic University of the Sacred Heart, via Musei 41, I-25121 Brescia, Italy}
\date{}
\maketitle
\begin{abstract}
In the Reduced Basis approximation of Stokes and Navier-Stokes problems, the Galerkin projection on the reduced spaces does not necessarily preserved the inf-sup stability even if the snapshots were generated through a stable full order method. Therefore, in this work we aim at building a stabilized Reduced Basis (RB) method for the approximation of unsteady Stokes and Navier-Stokes problems in parametric reduced order settings. This work extends the results presented for parametrized steady Stokes and Navier-Stokes problems in a work of ours \cite{Ali2018}. We apply classical residual-based stabilization techniques for finite element methods in full order, and then the RB method is introduced as Galerkin projection onto RB space. We compare this approach with supremizer enrichment options through several numerical experiments. We are interested to (numerically) guarantee the parametrized reduced inf-sup condition and to reduce the online computational costs.
\end{abstract}


\section{Introduction}
\label{sec:introduction}
In the finite element (FE) simulation of incompressible flows using a standard Galerkin formulation there are two possible sources of instabilities. One reason could be due to the presence of convection term which for high Reynolds number creates instability in numerical solution. Another source of instability could be due to the inappropriate choice of interpolating functions for velocity and pressure. Starting from early 70s, different researchers \cite{Roache1972,Christie1976,Heinrich1977,Heinrich1977(1),Hughes1978}  proposed several stabilized schemes to overcome stability issues. For instance, Hughes and Brooks \cite{Hughes1979, Hughes1980, Brooks1982} proposed to add artificial diffusion term acting only in the streamline direction and named this type of formulation as Streamline Upwind/Petrov Galerkin (SUPG) formulation. An extension of SUPG formulation is given by Hughes et al. \cite{Hughes1989} and is named as Galerkin Least Square (GALS) formulation. Later on Douglas-Wang \cite{Douglas1989} introduced the change of sign in GALS formulation.  A penalty method in which pressure is eliminated by penalizing the continuity equation and then retained in boundary condition was introduced by Hughes et al. \cite{Hughes1979(1)}. Hughes et al. \cite{Hughes1986} used equal order interpolation for velocity and pressure by perturbing the pressure test function with a gradient term to achieve the stability. A symmetric version of this method was given by Hughes and Franca \cite{HUGHES1987}. The SUPG method, first applied by Brooks and Hughes \cite{Brooks1982} to solve numerically the incompressible Navier-Stokes equations with high Reynolds number was later on extended by various researchers \cite{FRANCA1992(1),Douglas1989,Johnson1986,HANSBO1990,FRANCA1992,TEZDUYAR1992}.\par
Similarly the RB method for the Stokes \cite{Lov2006} and Navier-Stokes \cite{Lovgren2006} problems requires the fulfillment of discrete $\inf$-$\sup$ condition for reduced velocity and pressure spaces, respectively. In this paper we are not considering the convection dominated case, but we only focus on the $\inf$-$\sup$ stability at reduced order level. Previous works based on supremizer enrichment to cure the reduced $\inf$-$\sup$ condition are given by Rovas \cite{Rovas2003}, Rozza et al. \cite{Rozza2005a, Veroy2007}. Supremizer enrichment approach consists in the introduction of the inner pressure supremizer for the velocity-pressure stability of the RB spaces.
Several works on RB method for Stokes and Navier-Stokes problems using the pressure stabilization via the inner pressure supremizer operator are given by \cite{NegriManzoniRozza2015, Rozza2013,Simone2008,Simone2009,Manzoni2014,Veroy2005, Ballarin2015,BallarinChaconDelgadoGomezRozza2020,StabileBallari2018}.
\par
In our recent work on steady Stokes and Navier-Stokes problems  \cite{Ali2018, Ali2018b}, we proposed to use the classical residual based stabilization methods (such as SUPG, GALS and Douglas-Wang, mentioned above) to deal with the $\inf$-$\sup$ stability. This work is the continuation of proposed method to unsteady problems in parametric reduced order setting. We study the \textit{offline-online stabilization} \cite{PR014a} method, based on performing the Galerkin projections in both \textit{offline} and \textit{online} stage with respect to the consistent stabilized formulations, and the \textit{offline-only stabilization}, consisting in using the stabilized formulations only during the \textit{offline} stage and then projecting with respect to the standard formulation during the \textit{online} stage. We also guarantee the \textit{online} computational savings by reducing the dimension of the \textit{online} RB system, i.e, we show that with this approach it is possible to get the stable RB solution without the supremizer enrichment into velocity space.\par
This work has two parts: unsteady Stokes problem and unsteady Navier-Stokes problem. Further organization of this paper is as it follows: In section \ref{sec: unsteady_S_continuous_formulation} after recalling the unsteady Stokes problem, we present stabilized FE formulation and then its stabilized reduced basis (RB) formulation. Then we present some numerical results for unsteady Stokes problem in section \ref{eq:test1_unsteady} showing the error comparison between different stabilization and supremizer options.
\par
In section \ref{sec: Unsteady_NS} we follow a similar pattern for unsteady Navier-Stokes problem. We first define the full order FE formulation, followed by stabilized FE formulation, and then, we project onto RB space. Finally, we show some numerical results and discussions in section \ref{sec:Results_UNS}. The outcome of this work is summarized in section \ref{sec:Conclusion_NStokes}.

\section{Unsteady parametrized Stokes problem}
\label{sec: unsteady_S_continuous_formulation}
Let $\Omega\subset\mathbb{R}^2$, be a reference configuration, and we assume that current configuration $\Omega_o(\boldsymbol\mu)$ can be obtained as the image of map $\boldsymbol{T}(.;\boldsymbol\mu):\mathbb{R}^2\rightarrow\mathbb{R}^2,$ i.e. $\Omega_o(\boldsymbol\mu)=\boldsymbol{T}(\Omega;\boldsymbol\mu).$
The unsteady parametrized Stokes problem in current configuration reads as follows: find $\boldsymbol{u}_o(t;\boldsymbol\mu)\in \boldsymbol{V}$ and $p_o(t;\boldsymbol\mu)\in Q$ such that

\begin{equation}
\label{eq:unsteady_original}
\begin{cases}
\dfrac{\partial}{\partial{t}}\boldsymbol{u}_o-\nu\Delta{\boldsymbol{u}_o}+\nabla{p}_o= \boldsymbol{0} & \text{ in } \Omega_o(\boldsymbol\mu)\times\left(0,T\right),\  \\
\rm div\ {\boldsymbol{u}_o}=0 & \text{ in } \Omega_o(\boldsymbol\mu)\times\left(0,T\right),\\
\boldsymbol{u}_o=\boldsymbol{g}_D & \text{ on } \Gamma_{D,o}(\boldsymbol\mu)\times\left(0,T\right),\\
\boldsymbol{u}_o|_{t=0}=\boldsymbol{u_0} & \text{ on } \Gamma_{W,o}(\boldsymbol\mu),
\end{cases}
\end{equation}
where $(0,T)$ with $T>0$ is the time interval of interest, $\boldsymbol{u_0}\in L^2(\Omega)$ and $\nu$ is the viscosity of fluid. The boundary $\partial\Omega_o(\boldsymbol\mu)$ is divided into two parts in such a way that $\partial\Omega_o(\boldsymbol\mu)=\Gamma_{D,o}(\boldsymbol\mu)\cup\Gamma_{W,o}(\boldsymbol\mu),$ where $\Gamma_{D,o}(\boldsymbol\mu)$ is the Dirichlet boundary with non-homogeneous data and $\Gamma_{W,o}(\boldsymbol\mu)$ denotes the Dirichlet boundary with zero data.
\par We multiply \eqref{eq:unsteady_original} by velocity and pressure test functions $\boldsymbol{v}$ and $q$, respectively then integrating by parts, and tracing everything back onto the reference domain $\Omega,$ we obtain the following parametrized formulation of problem \eqref{eq:unsteady_original}:\par
for a given $\boldsymbol\mu\in\mathbb{P},$ find $\boldsymbol{u}(t;\boldsymbol\mu)\in\boldsymbol{V}$ and $p(t;\boldsymbol\mu)\in Q$ such that
\begin{equation}
\begin{cases}
m(\dfrac{\partial}{\partial{t}}\boldsymbol{u},\boldsymbol{v};\boldsymbol\mu)+a(\boldsymbol{u},\boldsymbol{v};\boldsymbol\mu)+b(\boldsymbol{v},p;\boldsymbol\mu)=F(\boldsymbol{v};\boldsymbol\mu) & \forall \, v\in \boldsymbol{V},t>0,\\
b(\boldsymbol{u},q;\boldsymbol\mu)=G(q;\boldsymbol\mu) & \forall \, q\in Q,t>0,\\
\boldsymbol{u}|_{t=0}=\boldsymbol{u_0}.
\end{cases}
\label{eq:unsteady_reference}
\end{equation}
We define the spaces $\boldsymbol{V} = L^2(\mathbb{R}^+;[H^1(\Omega)]^2) \cap C^0(\mathbb{R}^+;[L^2(\Omega)]^2)$ for velocity and $Q = L^2(\mathbb{R}^+;L_0^2(\Omega))$ for pressure on reference domain. Here, $H^1(\Omega)$ and $L^2(\Omega)$ are equipped with $H^1-$seminorm and $L^2-$norm respectively. Bilinear forms in \eqref{eq:unsteady_reference} are
\begin{equation}
\begin{split}
a(\boldsymbol{u},\boldsymbol{v};\boldsymbol\mu)&=\int_{\Omega}\dfrac{\partial\boldsymbol{u}}{\partial x_i}\kappa_{ij}(x;\boldsymbol\mu)\dfrac{\partial\boldsymbol{v}}{\partial x_j}d\boldsymbol{x}, \qquad
b(\boldsymbol{v},q;\boldsymbol\mu)=-\int_{\Omega}q \chi_{ij}(x;\boldsymbol\mu)\dfrac{\partial{v_j}}{\partial x_i}d\boldsymbol{x}, \\
m(\boldsymbol{u},\boldsymbol{v};\boldsymbol\mu)&=\int_{\Omega}\pi(\boldsymbol{x};\boldsymbol\mu)\boldsymbol{u}_i\boldsymbol{v}_id\boldsymbol{x}.
\end{split}
\label{eq: unsteady_bilinear_forms}
\end{equation}
We define the terms $F$ and $G$ in \eqref{eq:unsteady_reference} as:
\begin{equation}
\begin{split}
F(\boldsymbol{v};\boldsymbol\mu)&=-a(\boldsymbol{l(\mu)},\boldsymbol{v};\boldsymbol\mu),\\
G(q;\boldsymbol\mu)&=-b(\boldsymbol{l(\mu)},q;\boldsymbol\mu),
\end{split}
\label{eq: S_lift_functions}
\end{equation}
where we denote by $\boldsymbol{l}(\boldsymbol\mu)$ a parametrized lifting function such that  $\boldsymbol{l}(\boldsymbol\mu)|_{\Gamma_{D_g}}=\boldsymbol{g}_D(\boldsymbol\mu).$
\par
The tensors $\boldsymbol{\kappa}$, $\boldsymbol{\chi}$ and scalar $\pi$ encoding both physical and geometrical parametrization are defined as follows
\begin{equation}
\begin{split}
\boldsymbol{\kappa}(x;\boldsymbol\mu)&=\nu(J_T(x;\boldsymbol\mu))^{-1}(J_T(x;\boldsymbol\mu))^{-T}|J_T(X;\boldsymbol\mu)|, \\
\boldsymbol{\chi}(x;\boldsymbol\mu)&=(J_T(x;\boldsymbol\mu))^{-1}|J_T(X;\boldsymbol\mu)|, \\
\pi(\boldsymbol{x};\boldsymbol\mu)&=|J_T(X;\boldsymbol\mu)|,
\end{split}
\label{eq: unsteady_tensors}
\end{equation}
where $J_T\in \mathbb{R}^{2\times 2}$ is the Jacobian matrix of the map $\boldsymbol{T}(.;\boldsymbol\mu)$, and $|J_T|$ denotes the determinant.
\subsection{Semi-discrete Finite Element formulation}
\label{subsec: Algebraic_Stokes_unsteady_time}
The mixed Galerkin finite element semi-discretization \cite{Raviart1986,Max1989} of \eqref{eq:unsteady_reference} is defined as follows:\par
for a given $\boldsymbol\mu\in\mathbb{P},$ find $\boldsymbol{u}_h(t;\boldsymbol\mu)\in\boldsymbol{V}_h\subset\boldsymbol{V}$ and $p_h(t;\boldsymbol\mu)\in Q_h\subset Q$ such that
\begin{equation}
\begin{cases}
m(\dfrac{\partial}{\partial{t}}\boldsymbol{u}_h,\boldsymbol{v}_h;\boldsymbol\mu)+a(\boldsymbol{u}_h,\boldsymbol{v}_h;\boldsymbol\mu)+b(\boldsymbol{v}_h,p_h;\boldsymbol\mu)=F(\boldsymbol{v}_h;\boldsymbol\mu) & \forall \, \boldsymbol{v}_h\in \boldsymbol{V}_h,t>0,\\
b(\boldsymbol{u}_h,q_h;\boldsymbol\mu)=G(q_h;\boldsymbol\mu) & \forall \, q_h\in Q_h,t>0,\\
\boldsymbol{u}_{h}|_{t=0}=\boldsymbol{u}_{0,h}.
\end{cases}
\label{eq:unsteady_discrete}
\end{equation}
\par
We consider a partition of the interval $[0,T]$ into $K$ sub-intervals of equal length $\Delta{t}=T/K$ and $t^k=k\Delta{t}, 0\leq{k}\leq{K}.$ Applying the implicit Euler time discretization we obtain the following time discrete problem:\par
for a given $\boldsymbol\mu\in\mathbb{P},$ and $(\boldsymbol{u}_h^{k-1}(\boldsymbol\mu),p_h^{k-1}(\boldsymbol\mu)),$ find $\boldsymbol{u}_h^{k}(t;\boldsymbol\mu)\in\boldsymbol{V}_h$ and $p_h^{k}(t;\boldsymbol\mu)\in Q_h$ such that
\begin{equation}
\begin{cases}
\frac{1}{\Delta{t}}m(\boldsymbol{u}_h^{k},\boldsymbol{v}_h;\boldsymbol\mu)+a(\boldsymbol{u}_h^{k},\boldsymbol{v}_h;\boldsymbol\mu)+b(\boldsymbol{v}_h,p_h^{k};\boldsymbol\mu)=F(\boldsymbol{v}_h;\boldsymbol\mu)\\
+\frac{1}{\Delta{t}}m(\boldsymbol{u}_h^{k-1},\boldsymbol{v}_h;t^{k-1};\boldsymbol\mu) & \forall \, \boldsymbol{v}_h\in \boldsymbol{V}_h,\\
b(\boldsymbol{u}_h^{k},q_h;\boldsymbol\mu)=G(q_h;\boldsymbol\mu) & \forall \, q_h\in Q_h,\\
\boldsymbol{u}_{h}^{0}=\boldsymbol{u}_{0,h}.
\end{cases}
\label{eq:unsteady_discrete_time}
\end{equation}
\par
We provide the algebraic formulation of the semi-discrete problem \eqref{eq:unsteady_discrete}. The resulting ODE system is as follows:
\begin{equation}
\left[\begin{array}{cc} M(\boldsymbol\mu) & \boldsymbol{0}\\ \boldsymbol{0} & \boldsymbol{0} \end{array}\right]
\left[\begin{array}{c} \boldsymbol{\dot{U}}(t;\boldsymbol\mu) \\ \boldsymbol{\dot{P}}(t;\boldsymbol\mu) \end{array}\right] +
\left[\begin{array}{cc} A(\boldsymbol\mu) & B^T(\boldsymbol\mu)\\ B(\boldsymbol\mu) & \boldsymbol{0} \end{array}\right]
\left[\begin{array}{c} \boldsymbol{U}(t;\boldsymbol\mu) \\ \boldsymbol{P}(t;\boldsymbol\mu) \end{array}\right]
= \left[\begin{array}{c} \boldsymbol{\bar{f}}(\boldsymbol\mu) \\ \boldsymbol{\bar{g}}(\boldsymbol\mu) \end{array}\right]
\label{eq: Algebraic_Stokes_unsteady}
\end{equation}
for the vectors $\boldsymbol{U}=(u_h^{(1)},...,u_h^{({\mathcal{N}_u})})^T, \boldsymbol{P}=(p_h^{(1)},...,p_h^{({\mathcal{N}_p})})^T$, where for $1\leq i,j\leq \mathcal{N}_u$ and $1\leq k\leq \mathcal{N}_p.$ Let $\{\boldsymbol\phi_i^h\}_{i=1}^{\mathcal{N}_u}$ and $\{\psi_j^h\}_{j=1}^{\mathcal{N}_p}$ be basis functions of $\boldsymbol{V}_h$ and $Q_h$ respectively. We define the matrices
\begin{equation}
\begin{split}
(M(\boldsymbol\mu))_{ij}=m(\boldsymbol{\phi}_j^h,\boldsymbol{\phi}_i^h;\boldsymbol\mu),\qquad \left(A(\boldsymbol\mu)\right)_{ij}=a(\boldsymbol\phi_j^h,\boldsymbol\phi_i^h;\boldsymbol\mu),\\
\left(B(\boldsymbol\mu)\right)_{ki}=b(\boldsymbol\phi_i^h,\psi_k^h;\boldsymbol\mu),\qquad
(\boldsymbol{\bar{f}}(\boldsymbol\mu))_i=F(\boldsymbol\phi_i^h;\boldsymbol\mu), \\ (\boldsymbol{\bar{g}}(\boldsymbol\mu))_k=G(\psi_k^h;\boldsymbol\mu),
\end{split}
\label{eq: matrices_vectors}
\end{equation}

A key assumption for an efficient ROM evaluation is the capability to decouple the construction stage of the reduced order space (\textit{offline}) from evaluation stage (\textit{online}). We require that the matrices and vectors appearing in \eqref{eq: matrices_vectors} can be written as
\begin{equation}
\begin{split}
M(\boldsymbol\mu)=\sum_{q=1}^{Q_a}\Theta_{q}^{a}(\boldsymbol\mu)M^q,\quad
A(\boldsymbol\mu)=\sum_{q=1}^{Q_a}\Theta_{q}^{a}(\boldsymbol\mu)A^q, \qquad B(\boldsymbol\mu)=\sum_{q=1}^{Q_b}\Theta_{q}^{b}(\boldsymbol\mu)B^{q}, \\
\boldsymbol{\bar{f}}(\boldsymbol\mu)=\sum_{q=1}^{Q_f}\Theta_{q}^{f}(\boldsymbol\mu)\boldsymbol{\bar{f}}^q,\qquad  \boldsymbol{\bar{g}}(\boldsymbol\mu)=\sum_{q=1}^{Q_g}\Theta_{q}^{g}(\boldsymbol\mu)\boldsymbol{\bar{g}}^q.
\end{split}
\label{eq: affine}
\end{equation}
After applying the time discretization with implicit Euler scheme, the resulting algebraic formulation of \eqref{eq:unsteady_discrete_time} is
\begin{equation}
\begin{split}
\left[\begin{array}{cc} \dfrac{M(\boldsymbol\mu)}{\Delta t}+ A(\boldsymbol\mu)& B^T(\boldsymbol\mu)\\ B(\boldsymbol\mu) & \boldsymbol{0} \end{array}\right]
\left[\begin{array}{c} \boldsymbol{U}(t^{k};\boldsymbol\mu) \\ \boldsymbol{P}(t^{k};\boldsymbol\mu) \end{array}\right]
= \left[\begin{array}{c} \boldsymbol{\bar{f}}(\boldsymbol\mu) \\ \boldsymbol{\bar{g}}(\boldsymbol\mu) \end{array}\right] \\+
\left[\begin{array}{cc} \dfrac{M(\boldsymbol\mu)}{\Delta t} & \boldsymbol{0}\\ \boldsymbol{0} & \boldsymbol{0} \end{array}\right]
\left[\begin{array}{c} \boldsymbol{U}(t^{k-1};\boldsymbol\mu) \\ \boldsymbol{P}(t^{k-1};\boldsymbol\mu) \end{array}\right].
\end{split}
\label{eq: Algebraic_Stokes_unsteady_time}
\end{equation}
For a stable solution the FE spaces $\boldsymbol{V}_h$ and $Q_h$  have to fulfill the following parametrized inf-sup stability condition (LBB) \cite{QV}:
\begin{equation}
\exists \beta_{0}(\boldsymbol\mu)>0: \beta_{h}(\boldsymbol{\mu})= \inf_{q_h\in Q_{h}}\sup_{\boldsymbol{v}_h\in \boldsymbol{V}_{h}}\frac{b(\boldsymbol{v}_h,q_h;\boldsymbol\mu)}{\|\boldsymbol{v}_h\|_{\boldsymbol{V}_h}\|q_h\|_{Q_h}}\geq \beta_{0}(\boldsymbol\mu) \quad\forall \boldsymbol\mu\in \mathbb{P}.
\label{eq: inf-sup}
\end{equation}
This relation holds if, e.g., the Taylor-Hood $(\mathbb{P}_2/\mathbb{P}_1)$ FE spaces are chosen. It is important to mention that condition \eqref{eq: inf-sup} does not hold in case of equal order FE spaces  $(\mathbb{P}_k/\mathbb{P}_k), k\geq 1$ and for lowest order element $(\mathbb{P}_1/\mathbb{P}_0)$. Therefore, in such situations we need to introduce some additional stabilization terms, as in the following.
\subsection{Stabilized Finite Element formulation}
\label{subsec: Algebraic_Stokes_unsteady_stabFE}
Let us modify equation \eqref{eq:unsteady_discrete} by adding the stabilization terms. We read the modified formulation as follows: for a given $\boldsymbol\mu\in\mathbb{P},$ find $\boldsymbol{u}_h(t;\boldsymbol\mu)\in\boldsymbol{V}_h$ and $p_h(t;\boldsymbol\mu)\in Q_h$ such that
\begin{equation}
\begin{cases}
m(\dfrac{\partial}{\partial{t}}\boldsymbol{u}_h,\boldsymbol{v}_h;\boldsymbol\mu)+a(\boldsymbol{u}_h,\boldsymbol{v}_h;\boldsymbol\mu)+b(\boldsymbol{v}_h,p_h;\boldsymbol\mu)=F(\boldsymbol{v}_h;\boldsymbol\mu) & \forall \, \boldsymbol{v}_h\in \boldsymbol{V}_h,t>0,\\
b(\boldsymbol{u}_h,q_h;\boldsymbol\mu)-s^{ut,q}_{h}(\boldsymbol{u}_{h},q_{h};\boldsymbol\mu)-s^{u,q}_{h}(\boldsymbol{u}_{h},q_{h};\boldsymbol\mu)-s^{p,q}_{h}(p_h,q_{h};\boldsymbol\mu)= G(q_h;\boldsymbol\mu) & \forall \, q_h\in Q_h,t>0,\\
\boldsymbol{u}_{h}|_{t=0}=\boldsymbol{u}_{0,h},
\end{cases}
\label{eq:unsteady_discrete_stab}
\end{equation}
where $s^{ut,q}_{h}(.,.;\boldsymbol\mu), s^{u,q}_{h}(.,.;\boldsymbol\mu)$ and $s^{p,q}_{h}(.,.;\boldsymbol\mu)$ are the stabilization terms. For a detail discussion on the choice of stabilization terms, we refer to recent work of ours \cite{Ali2018} and references therein. In this case we prefer to chose the stabilization technique given by Hughes et al. \cite{Hughes1986}:
\begin{equation}
s^{ut,q}_{h}(\boldsymbol{u}_{h},q_{h};\boldsymbol\mu):=\delta\sum_{K}h_K^{2}\int_K(\dfrac{\partial}{\partial{t}}\boldsymbol{u}_h,\nabla{q}_{h}),
\label{eq:our_choice_unsteady}
\end{equation}
\begin{equation}
s^{u,q}_{h}(\boldsymbol{u}_{h},q_{h};\boldsymbol\mu):=\delta\sum_{K}h_K^{2}\int_K(-\nu\Delta\boldsymbol{u}_{h},\nabla{q}_{h}),
\label{eq:our_choice_unsteady_1}
\end{equation}
and
\begin{equation}
s^{p,q}_{h}(p_h,q_{h};\boldsymbol\mu):=\delta\sum_{K}h_K^{2}\int_K(\nabla{p}_{h},\nabla{q}_{h}),
\label{eq:our_choice_unsteady_2}
\end{equation}
Therefore, the stabilized algebraic system can be written as:
\begin{equation}
\left[\begin{array}{cc} M(\boldsymbol\mu) & \boldsymbol{0}\\ \tilde{M}(\boldsymbol\mu) & \boldsymbol{0} \end{array}\right]
\left[\begin{array}{c} \boldsymbol{\dot{U}}(t;\boldsymbol\mu) \\ \boldsymbol{\dot{P}}(t;\boldsymbol\mu) \end{array}\right] +
\left[\begin{array}{cc} A(\boldsymbol\mu) & B^T(\boldsymbol\mu)\\ \tilde{B}(\boldsymbol\mu) & -S(\boldsymbol\mu) \end{array}\right]
\left[\begin{array}{c} \boldsymbol{U}(t;\boldsymbol\mu) \\ \boldsymbol{P}(t;\boldsymbol\mu) \end{array}\right]
= \left[\begin{array}{c} \boldsymbol{\bar{f}}(\boldsymbol\mu) \\ \boldsymbol{\bar{g}}(\boldsymbol\mu) \end{array}\right],
\label{eq: Algebraic_Stokes_unsteady_stab}
\end{equation}
where $\tilde{M}(\boldsymbol\mu), \tilde{B}(\boldsymbol\mu)$ and $-S(\boldsymbol\mu)$ contains the stabilization effects \cite{Hughes1986}, and defined as follows:
\begin{equation}
\begin{split}
\left(\tilde{M}(\boldsymbol\mu)\right)_{ki}&=s^{ut,q}_{h}(\boldsymbol\phi_i^h,\psi_k^h;\boldsymbol\mu),\qquad
\left(\tilde{B}(\boldsymbol\mu)\right)_{ki}=b(\boldsymbol\phi_i^h,\psi_k^h;\boldsymbol\mu)+s^{u,q}_{h}(\boldsymbol\phi_i^h,\psi_k^h;\boldsymbol\mu),\\
(S(\boldsymbol\mu))_{ij}&=s^{p,q}_{h}(\psi^{h}_j,\psi^{h}_i;\boldsymbol\mu),\quad \text{ for }\quad 1\leq i,j\leq\mathcal{N}_u, 1\leq k\leq\mathcal{N}_p,
\end{split}
\label{eq: Stokes_stab_matrices}
\end{equation}
After applying the time discretization with implicit Euler scheme, the system \eqref{eq: Algebraic_Stokes_unsteady_stab} becomes
\begin{equation}
\begin{split}
\left[\begin{array}{cc} \dfrac{M(\boldsymbol\mu)}{\Delta t}+ A(\boldsymbol\mu)& B^T(\boldsymbol\mu)\\ \tilde{B}(\boldsymbol\mu)+\dfrac{\tilde{M}(\boldsymbol\mu)}{\Delta t} & -S(\boldsymbol\mu) \end{array}\right]
\left[\begin{array}{c} \boldsymbol{U}(t^{k};\boldsymbol\mu) \\ \boldsymbol{P}(t^{k};\boldsymbol\mu) \end{array}\right]
= \left[\begin{array}{c} \boldsymbol{\bar{f}}(\boldsymbol\mu) \\ \boldsymbol{\bar{g}}(\boldsymbol\mu) \end{array}\right] \\+
\left[\begin{array}{cc} \dfrac{M(\boldsymbol\mu)}{\Delta t} & \boldsymbol{0}\\ \dfrac{\tilde{M}(\boldsymbol\mu)}{\Delta t} & \boldsymbol{0} \end{array}\right]
\left[\begin{array}{c} \boldsymbol{U}(t^{k-1};\boldsymbol\mu) \\ \boldsymbol{P}(t^{k-1};\boldsymbol\mu) \end{array}\right].
\end{split}
\label{eq: Algebraic_Stokes_unsteady_stabFE}
\end{equation}
The stabilized formulation requires the FE spaces to fulfill the following modified $\inf$-$\sup$ condition \cite{boffi2013mixed,Burman2009,Becker2001} after adding some additional stabilization terms:
\begin{equation}
\exists \beta_{0}(\boldsymbol{\mu})>0:  \sup_{{\boldsymbol{v}_h}\in \boldsymbol{V}_{h}}\frac{b(\boldsymbol{v}_h,q_h;\boldsymbol\mu)}{\|\nabla{\boldsymbol{v}_h}\|}+s^{p,q}_{h}\left(q_{h},q_{h};\boldsymbol\mu\right)^{1/2}\geq \beta_0(\boldsymbol\mu)\|q_h\|, \forall q_h\in Q_{h}.
\label{eq:modified_inf-sup}
\end{equation}
\subsection{Reduced Basis formulation}
In this section we present the RB formulation of the unsteady Stokes problem formulated in section \ref{subsec: Algebraic_Stokes_unsteady_time}. Let us define the parameter sample $s_{N}=\left\lbrace\boldsymbol\mu^1,...,\boldsymbol\mu^N\right\rbrace$, where $\boldsymbol\mu^n\in \mathbb{P}$. The reduced basis approximation is based on an $N$-dimensional reduced basis spaces $\boldsymbol{V}_{N}$ and $Q_{N}$ generated by a sampling procedure which combines spatial snapshots in time and parameter space in an optimal way. In particular, in our case we have used the POD-Greedy algorithm \cite{Haasdonk2013} for snapshots selection to generate the reduced spaces. Reduced basis velocity and pressure spaces are
\begin{equation}
\boldsymbol{V}_N = \text{ span} \left\lbrace\text{POD}(\boldsymbol{u}_{h}(t^{k};\boldsymbol\mu^n)), 1\leq k\leq K,1\leq n\leq N_u\right\rbrace,
\label{eq:velocity}
\end{equation}
\begin{equation}
Q_N=\text{ span} \left\lbrace\text{POD}(p_{h}(t^{k};\boldsymbol\mu^n)), 1\leq k\leq K,1\leq n\leq N_p\right\rbrace.
\label{eq:pressure_t}
\end{equation}
We introduce the supremizer operator $T^{\boldsymbol\mu}:Q_h\rightarrow\boldsymbol{V}_h$ defined as follows:
\begin{equation}
(T^{\boldsymbol\mu}q_h,\boldsymbol{v}_h)_{\boldsymbol{V}}=b(\boldsymbol{v}_h,q_h;\boldsymbol\mu), \quad \forall \boldsymbol{v}\in \boldsymbol{V}_h.
\label{eq: supremizer}
\end{equation}
which is evaluated for $\boldsymbol\mu = \boldsymbol\mu^n$ and the corresponding pressure snapshot $q_h^{k} := p_{h}(t^{k};\boldsymbol\mu^n)$, $n = 1, \hdots, N$, to obtain $N$ supremizer snapshots.
Afterwards, the RB velocity space $\boldsymbol{V}_N$ is enriched with the supremizer snapshots. We denote the enriched RB velocity space by $\tilde{\boldsymbol{V}}_N$, defined as:
\begin{equation}
\tilde{\boldsymbol{V}}_N = \text{span}\left\lbrace\text{POD}(\boldsymbol{u}_{h}(t^{k};\boldsymbol\mu^n)), 1\leq n\leq N_u;\text{POD}(T^{\boldsymbol\mu}q_{h}(t^{k};\boldsymbol\mu^n)), 1\leq n\leq N_s\right\rbrace, 
\label{eq: sup_velocity}
\end{equation}
where $N_s\leq N_p$ denotes the number of supremizer snapshots. Now the reduced basis formulation corresponding to semi-discrete FE formulation \eqref{eq:unsteady_discrete} can be written as: for any $\boldsymbol\mu\in\mathbb{P},$ find $\boldsymbol{u}_N(t;\boldsymbol\mu)\in\boldsymbol{V}_N$ and $p_N(t;\boldsymbol\mu)\in Q_N$ such that
\begin{equation}
\begin{cases}
m(\dfrac{\partial}{\partial{t}}\boldsymbol{u}_N,\boldsymbol{v}_N;\boldsymbol\mu)+a(\boldsymbol{u}_N,\boldsymbol{v}_N;\boldsymbol\mu)+b(\boldsymbol{v}_N,p_N;\boldsymbol\mu)=F(\boldsymbol{v}_N;\boldsymbol\mu) & \forall \, \boldsymbol{v}_N\in \boldsymbol{V}_N,\\
b(\boldsymbol{u}_N,q_N;\boldsymbol\mu)=G(q_N;\boldsymbol\mu) & \forall \, q_N\in Q_N,\\
\boldsymbol{u}_{N}|_{t=0}=\boldsymbol{u}_{0,N}.
\end{cases}
\label{eq:unsteady_discrete_RB}
\end{equation}
In the online stage, the algebraic formulation of resulting reduced order approximation for any $\boldsymbol\mu\in\mathbb{P}$ is given by
\begin{equation}
\begin{split}
\left[\begin{array}{cc} \dfrac{M_N(\boldsymbol\mu)}{\Delta t}+ A_N(\boldsymbol\mu)& B_N^T(\boldsymbol\mu)\\ B_N(\boldsymbol\mu) & \boldsymbol{0} \end{array}\right]
\left[\begin{array}{c} \boldsymbol{U}_N(t^{k};\boldsymbol\mu) \\ \boldsymbol{P}_N(t^{k};\boldsymbol\mu) \end{array}\right]
= \left[\begin{array}{c} \boldsymbol{\bar{f}}_N(\boldsymbol\mu) \\ \boldsymbol{\bar{g}}_N(\boldsymbol\mu) \end{array}\right] \\+
\left[\begin{array}{cc} \dfrac{M_N(\boldsymbol\mu)}{\Delta t} & \boldsymbol{0}\\ \boldsymbol{0} & \boldsymbol{0} \end{array}\right]
\left[\begin{array}{c} \boldsymbol{U}_N(t^{k-1};\boldsymbol\mu) \\ \boldsymbol{P}_N(t^{k-1};\boldsymbol\mu) \end{array}\right],
\end{split}
\label{eq: Algebraic_Stokes_unsteady_time_RB}
\end{equation}
where the reduced order matrices are defined as:
\begin{equation}
\begin{split}
M_N(t;\boldsymbol\mu)=Z^T_{u,s}M(t;\boldsymbol\mu)Z_{u,s}, \quad  A_N(\boldsymbol\mu)=Z^T_{u,s}A(\boldsymbol\mu)Z_{u,s}, \quad B_N(\boldsymbol\mu)=Z^T_{p}B(\boldsymbol\mu)Z_{u,s}, \\ \boldsymbol{\bar{f}}_N(\boldsymbol\mu)=Z^T_{u,s}\boldsymbol{\bar{f}}(\boldsymbol\mu),\quad
\boldsymbol{\bar{g}}_N(\boldsymbol\mu)=Z^T_{p}\boldsymbol{\bar{g}}(\boldsymbol\mu),
\end{split}
\end{equation}
with $Z_{u,s}$ being the velocity snapshot matrix including the supremizer solutions, $Z_{p}$ denotes the pressure snapshot matrix.
Moreover, thanks to the affine parametric dependence \eqref{eq: affine}, we need to store only the matrices and vectors
\begin{equation}
A^q_N=Z^T_{u,s}A^qZ_{u,s}, \quad B^q_N=Z^T_{p}B^qZ_{u,s}, \quad \boldsymbol{\bar{f}}^q_N=Z^T_{u,s}\boldsymbol{\bar{f}}^q, \quad
\boldsymbol{\bar{g}}^q_N=Z^T_{p}\boldsymbol{\bar{g}}^q.
\end{equation}
The store data structures do not depend explicitly on time because the temporal dependence is stored in the multiplicative factors $\Theta(t;\boldsymbol\mu).$ Therefore, $A^q,B^q,\boldsymbol{\bar{f}}^q,\boldsymbol{\bar{g}}^q$ are independent of both $\boldsymbol\mu$ and $t$.
\label{subsec: unsteadyS_RB_formulation}
\subsection{Stabilized Reduced Basis formulation}
In this section we present the stabilized RB model for unsteady Stokes problem derived from the stabilized FE problem \eqref{eq:unsteady_discrete_stab}. The stabilized RB approximation of velocity and pressure field obtained by means of Galerkin projection on reduced spaces reads:\par
for any $\boldsymbol\mu\in\mathbb{P},$ find $\boldsymbol{u}_N(t;\boldsymbol\mu)\in\boldsymbol{V}_N$ and $p_N(t;\boldsymbol\mu)\in Q_N$ such that
\begin{equation}
\begin{cases}
m(\dfrac{\partial}{\partial{t}}\boldsymbol{u}_N,\boldsymbol{v}_N;\boldsymbol\mu)+a(\boldsymbol{u}_N,\boldsymbol{v}_N;\boldsymbol\mu)+b(\boldsymbol{v}_N,p_N;\boldsymbol\mu)=F(\boldsymbol{v}_N;\boldsymbol\mu) & \forall \, \boldsymbol{v}_N\in \boldsymbol{V}_N,t>0,\\
b(\boldsymbol{u}_N,q_N;\boldsymbol\mu)-s^{ut,q}_{N}(\boldsymbol{u}_{N},q_{N};\boldsymbol\mu)-s^{u,q}_{N}(\boldsymbol{u}_{N},q_{N};\boldsymbol\mu)-s^{p,q}_{N}(p_N,q_N;\boldsymbol\mu)= G(q_N;\boldsymbol\mu) & \forall \, q_N\in Q_N,t>0,\\
\boldsymbol{u}_{N}|_{t=0}=\boldsymbol{u}_{0,N},
\end{cases}
\label{eq:unsteady_discrete_stab_RB}
\end{equation}
where $s^{ut,q}_{N}(.,.;\boldsymbol\mu), s^{u,q}_{N}(.,.;\boldsymbol\mu)$ and $s^{p,q}_{N}(.,.;\boldsymbol\mu)$ are the reduced order stabilization terms defined as:
\begin{equation}
s^{ut,q}_{N}(\boldsymbol{u}_{N},q_{N};\boldsymbol\mu):=\delta\sum_{K}h_K^{2}\int_K(\dfrac{\partial}{\partial{t}}\boldsymbol{u}_N,\nabla{q}_{N}),
\label{eq:our_choice_unsteady_RB}
\end{equation}
\begin{equation}
s^{u,q}_{N}(\boldsymbol{u}_{N},q_{N};\boldsymbol\mu):=\delta\sum_{K}h_K^{2}\int_K(-\nu\Delta\boldsymbol{u}_{N},\nabla{q}_{N}),
\label{eq:our_choice_unsteady_1_RB}
\end{equation}
and
\begin{equation}
s^{p,q}_{N}(p_N,q_N;\boldsymbol\mu):=\delta\sum_{K}h_K^{2}\int_K(\nabla{p}_{N},\nabla{q}_{N}),
\label{eq:our_choice_unsteady_2_RB}
\end{equation}
Finally, we write the reduced order stabilized formulation of unsteady FE stabilized Stokes problem \eqref{eq: Algebraic_Stokes_unsteady_stabFE} in compact form as:
\begin{equation}
\begin{split}
\left[\begin{array}{cc} \dfrac{M_N(\boldsymbol\mu)}{\Delta t}+ A_N(\boldsymbol\mu)& B_N^T(\boldsymbol\mu)\\ \tilde{B}_N(\boldsymbol\mu)+\dfrac{\tilde{M}_N(\boldsymbol\mu)}{\Delta t} & -S_N(\boldsymbol\mu) \end{array}\right]
\left[\begin{array}{c} \boldsymbol{U}_N(t^{k};\boldsymbol\mu) \\ \boldsymbol{P}_N(t^{k};\boldsymbol\mu) \end{array}\right]
= \left[\begin{array}{c} \boldsymbol{\bar{f}}_N(\boldsymbol\mu) \\ \boldsymbol{\bar{g}}_N(\boldsymbol\mu) \end{array}\right] \\+
\left[\begin{array}{cc} \dfrac{M_N(\boldsymbol\mu)}{\Delta t} & \boldsymbol{0}\\ \dfrac{\tilde{M}_N(\boldsymbol\mu)}{\Delta t} & \boldsymbol{0} \end{array}\right]
\left[\begin{array}{c} \boldsymbol{U}_N(t^{k-1};\boldsymbol\mu) \\ \boldsymbol{P}_N(t^{k-1};\boldsymbol\mu) \end{array}\right].
\end{split}
\label{eq: Algebraic_Stokes_unsteady_stabRB}
\end{equation}
where $\tilde{M}_N(\boldsymbol\mu),$ $\tilde{B}_N(\boldsymbol\mu)$ and $S_N(\boldsymbol\mu)$ are RB stabilization matrices defined as:
\begin{equation}
\tilde{M}_N(\boldsymbol\mu)=Z^T_{p}\tilde{M}(\boldsymbol\mu)Z_{u,s}, \quad \tilde{B}_N(\boldsymbol\mu)=Z^T_{p}\tilde{B}(\boldsymbol\mu)Z_{u,s}, \quad S_N(\boldsymbol\mu)=Z^T_{p}S(\boldsymbol\mu)Z_{p},
\label{eq:RBmatrices_SteadyStokes_STAB}
\end{equation}
We also define the reduced order generalized $\inf$-$\sup$ condition
\begin{equation}
\exists \beta_{0,N}(\boldsymbol{\mu})>0: \sup_{{\boldsymbol{v}_N}\in \boldsymbol{V}_{N}}\frac{b(\boldsymbol{v}_N,q_N;\boldsymbol\mu)}{\|\nabla{\boldsymbol{v}_N}\|}+s^{p,q}_{N}(q_{N},q_{N};\boldsymbol\mu)^{1/2}\geq \beta_{0,N}(\boldsymbol{\mu}) \|q_N\|, \forall q_N\in Q_{N},
\label{eq:modified_inf-sup_RB}
\end{equation}
where $s^{p,q}_{N}(.,.;\boldsymbol\mu)$ is due to the addition of stabilization terms in RB formulation.
\par
For a detailed discussion about the combination of supremizer and stabilization approaches to fulfill the reduced inf-sup condition \eqref{eq:modified_inf-sup_RB}, we refer to our recent work \cite{Ali2018}. Here, we discuss and compare the following options using unstable FE pair $\mathbb{P}_k/\mathbb{P}_k$ :
\begin{itemize}
\item for \textit{offline-online stabilization} with supremizer we solve the stabilized system \eqref{eq: Algebraic_Stokes_unsteady_stabFE} in the \textit{offline} stage and stabilized RB system \eqref{eq: Algebraic_Stokes_unsteady_stabRB} in the \textit{online} stage; and the velocity space in this case is enriched with supremizer solutions, given by \eqref{eq: sup_velocity};
\item for \textit{offline-online stabilization} without supremizer we solve the stabilized system \eqref{eq: Algebraic_Stokes_unsteady_stabFE} in the \textit{offline} stage and stabilized RB system \eqref{eq: Algebraic_Stokes_unsteady_stabRB} in the \textit{online} stage; but the velocity space in this case is given by \eqref{eq:velocity};
\item for \textit{offline-only stabilization} with supremizer we solve the stabilized system \eqref{eq: Algebraic_Stokes_unsteady_stabFE} in the \textit{offline} stage and non-stabilized RB system \eqref{eq: Algebraic_Stokes_unsteady_time_RB} in the \textit{online} stage; and the velocity space in this case is enriched with supremizer solutions, given by \eqref{eq: sup_velocity};
\end{itemize}
\label{subsec:RBStokes_unsteady_stab}

\section{Numerical results and discussion}
\label{eq:test1_unsteady}
In this section, we present several numerical results for stabilized reduced order model for unsteady Stokes problem developed in section \ref{sec: unsteady_S_continuous_formulation}.
\par
\begin{figure}[H]
\centering
\includegraphics[width=0.45\textwidth]{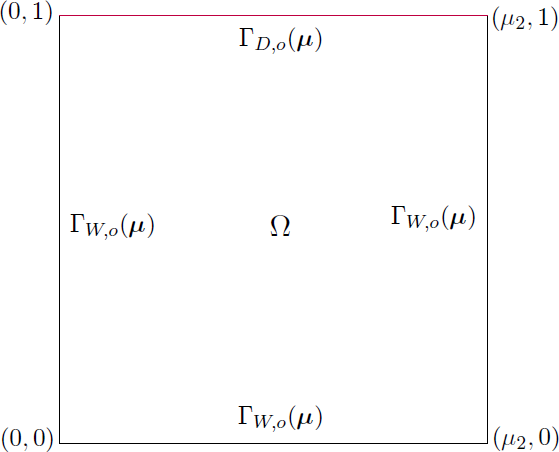}
\caption{Parametrized domain}
\label{fig:Domain}
\end{figure}

We set the parametrized domain $\Omega_o(\boldsymbol\mu)=(0,\mu_2)\times(0,1)$, where we define $\boldsymbol{\mu}=(\mu_1,\mu_2)$ such that $\mu_1$ is physical parameter (kinematic viscosity of fluid) and $\mu_2$ is geometrical parameter (length of domain). Main goal is to see the effect of geometrical parameter on the velocity and pressure. Parametrized domain is shown in Fig. \ref{fig:Domain}.
\par
We consider a partition of the boundary $\partial\Omega$ into $\Gamma_{D,o}(\boldsymbol\mu)\cup\Gamma_{W,o}(\boldsymbol\mu)$, where we have the homogeneous Dirichlet condition on $\Gamma_{W,o}(\boldsymbol\mu)$ and non-homogeneous Dirichlet condition on $\Gamma_{D,o}(\boldsymbol\mu)$.
\par
\subsection{Numerical results for $\mathbb{P}_k/\mathbb{P}_k$ (for $k=1, 2$)}
\label{subsec: PkPk}
The aim of present subsection is to show and discuss some numerical results for unsteady parametrized Stokes problem \eqref{eq:test1_unsteady} using Franca-Hughes stabilization \cite{Hughes1986}.
\par
In Table \ref{tab:comput_unsteady} we show the details of parameter ranges in \textit{offline}, \textit{online} stages; and other information about the \textit{offline} stage.\par
\begin{figure}[H]
\begin{tabular}{ll}
\includegraphics[width = 0.5\textwidth, keepaspectratio]{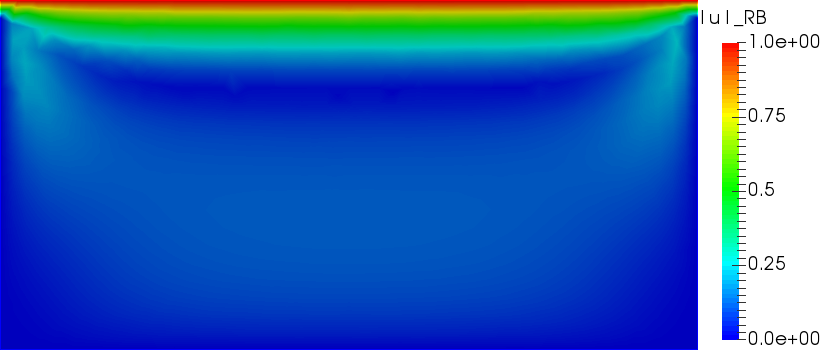} & \includegraphics[width = 0.5\textwidth, keepaspectratio]{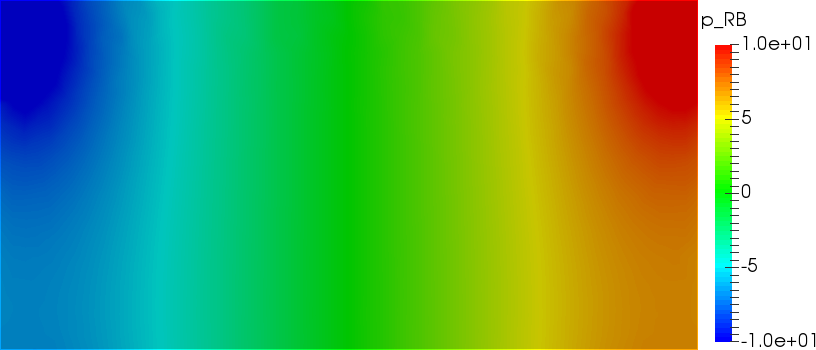}\\
 & \\
\includegraphics[width = 0.5\textwidth, keepaspectratio]{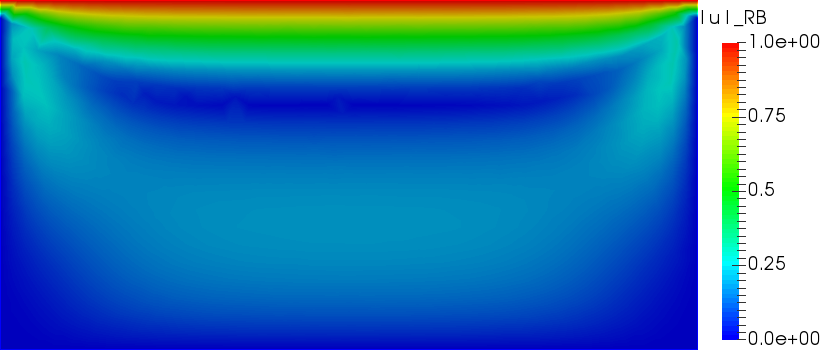} & \includegraphics[width = 0.5\textwidth, keepaspectratio]{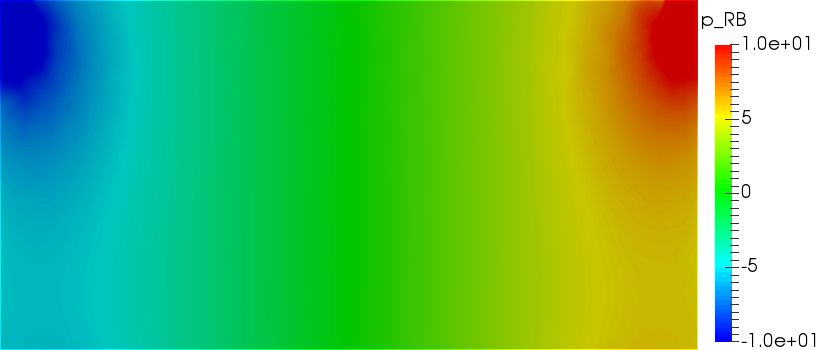}\\
 & \\
   \includegraphics[width = 0.5\textwidth, keepaspectratio]{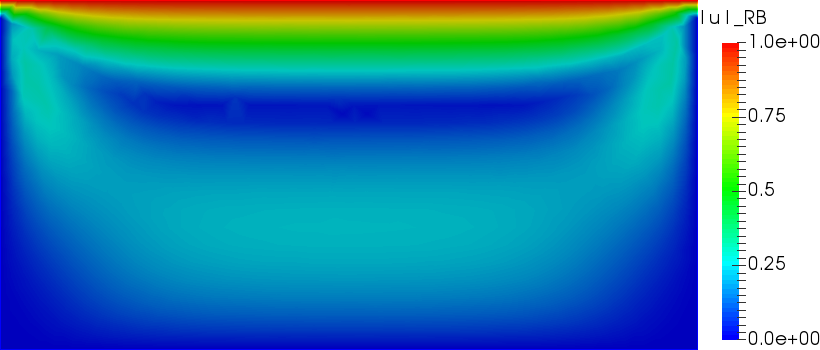} & \includegraphics[width = 0.5\textwidth, keepaspectratio]{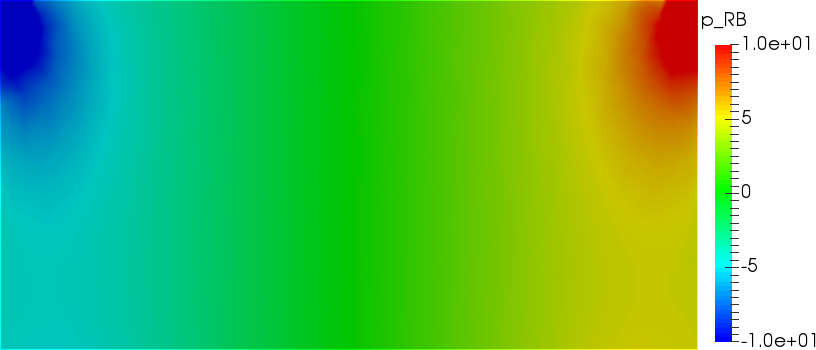}\\
  & \\
    \includegraphics[width = 0.5\textwidth, keepaspectratio]{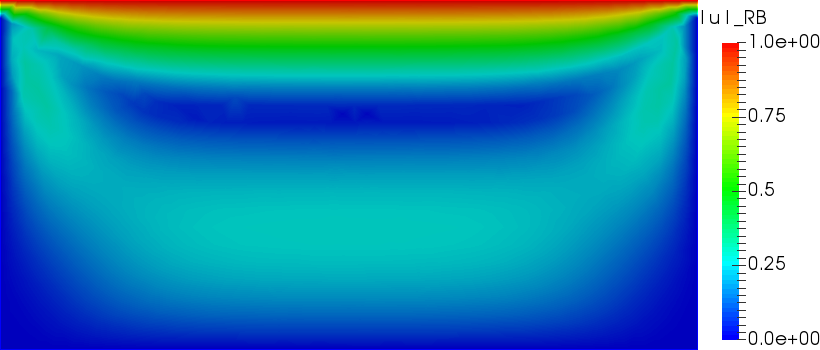} & \includegraphics[width = 0.5\textwidth, keepaspectratio]{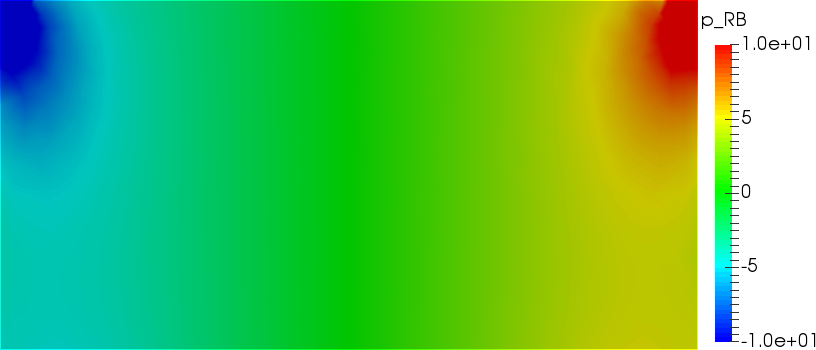}\\
  & \\
    \includegraphics[width = 0.5\textwidth, keepaspectratio]{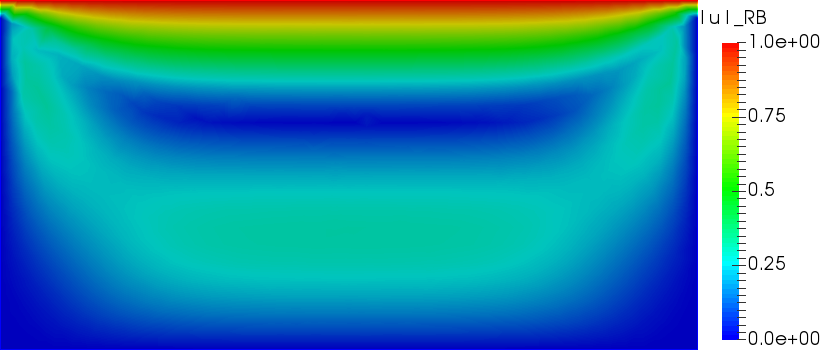} & \includegraphics[width = 0.5\textwidth, keepaspectratio]{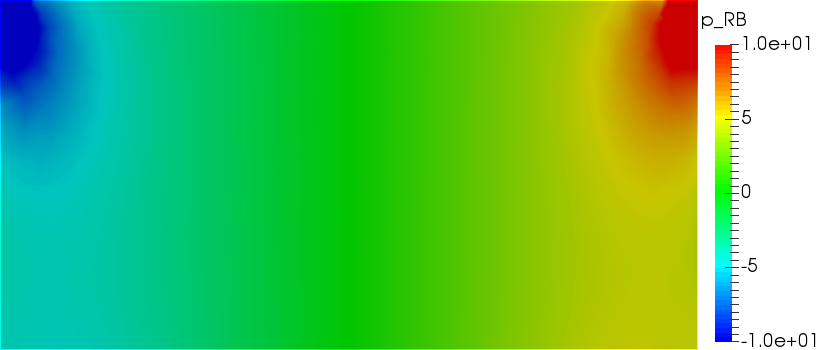}\\
\end{tabular}
 \caption{Stokes problem: Franca-Hughes stabilization with $\mathbb{P}_2/\mathbb{P}_2$ FE pair; RB solutions for Velocity field (left) and Pressure field (right) at different time step from top to bottom; $t=0.02,0.04,0.06,0.1,0.12$, $N_u=N_p=30.$}
\label{fig:unsteady_stokes Sol}
\end{figure}

In Fig. \ref{fig:unsteady_stokes Sol} we show the RB solutions for velocity and pressure at different time steps using the \textit{offline-online stabilization} without supremizer. We observe that as the time increases, both velocity and pressure fields are converging to steady state solutions. We have similar results with \textit{offline-online stabilization} with supremizer that we do not show here. \par
\begin{figure}[H]
  \centering
  \begin{subfigure}{0.48\textwidth}
    \centering\includegraphics[width=\textwidth]{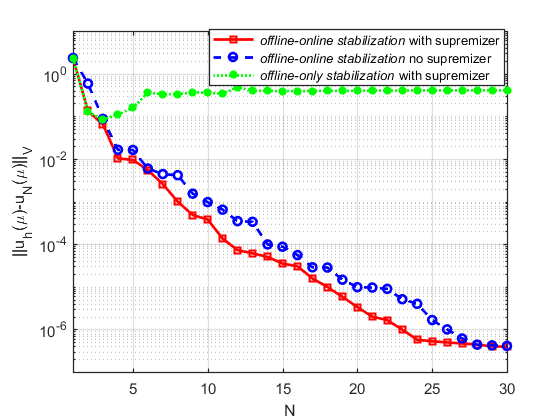}
  \end{subfigure}%
  \quad%
  \begin{subfigure}{0.48\textwidth}
    \centering\includegraphics[width=\textwidth]{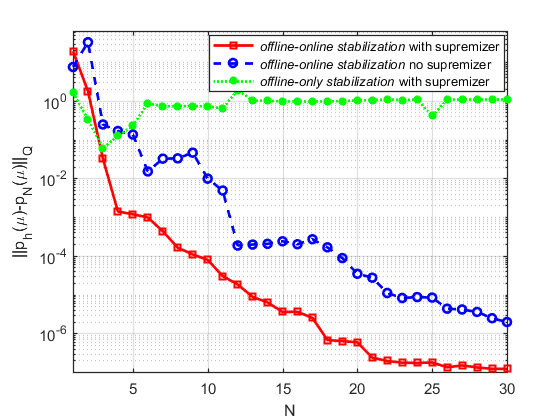}
  \end{subfigure}
  \caption{Stokes problem: Franca-Hughes stabilization with $\mathbb{P}_2/\mathbb{P}_2$ on cavity flow; $L^2$-error in time for velocity (left) and pressure (right) with stabilization coefficient $\delta=0.05$ and $\Delta{t}=0.02.$}
    \label{fig: UL2_alpha0_05}
\end{figure}

Figure \ref{fig: UL2_alpha0_05} shows the error between FE and RB solutions for velocity (left) and pressure (right), respectively. From these plots we observe that the \textit{offline-online stabilization} with and without supremizer show the same convergence behavior in case of velocity but in case of pressure, supremizer is improving the \textit{offline-online stabilization} up to one order of magnitude. We have similar results for $\mathbb{P}_1/\mathbb{P}_1$ FE pair. This property will be much important in case of coupling conditions in multi-physics involving pressure, for example, since we may guarantee a better accuracy. In contrast, the \textit{offline-only stabilization} with supremizer option has poor performance for both velocity and pressure. From table \ref{tab:comput_unsteady} we see that the computation time of \textit{offline-online stabilization} without supremizer is less than the computation time of \textit{offline-online stabilization} with supremizer in both \textit{offline} and \textit{online stages}.

\begin{center}
 \begin{tabular}{||c | c||}
 \hline\hline
 Number of Parameters & 2: $\mu_1$(viscosity), $\mu_2$(domain's length) \\
 \hline
 $\mu_1$ range \textit{offline} & [0.25,0.75]  \\
 \hline
 $\mu_2$ range \textit{offline} & [1,2]  \\
 \hline
$\mu_1$ value \textit{online} & 0.57  \\
 \hline
$\mu_2$ value \textit{online} & 1.78  \\
 \hline
Final time & 0.2\\
 \hline
Time step $\Delta t$ & 0.02\\
 \hline
$N_{train}$ & 25\\
\hline
$N_{max}$ & 25\\
\hline
Stabilization coefficient $\delta$ & 0.05\\
\hline
\multirow{2}{*}{FE degrees of freedom} & $6222$ ($\mathbb{P}_1/{\mathbb{P}_1}$)\\
& $18300$ ($\mathbb{P}_2/{\mathbb{P}_2}$)\\
 \hline
RB dimension & $N_u=N_s=N_p=30$\\
 \hline
Computation time ($\mathbb{P}_2/{\mathbb{P}_1}$) & $1780s$ (\textit{offline}), $300s$ (\textit{online}) with supremizer\\
 \hline
\multirow{3}{*}{\textit{Offline} time ($\mathbb{P}_1/{\mathbb{P}_1}$)} & $1046s$ (\textit{offline-online stabilization} with supremizer)\\
& $738s$ (\textit{offline-online stabilization} without supremizer)\\
& $980s$ (\textit{offline-only stabilization} with supremizer)\\
 \hline
\multirow{3}{*}{\textit{Offline} time ($\mathbb{P}_2/{\mathbb{P}_2}$)} & $2260s$ (\textit{offline-online stabilization} with supremizer)\\
& $1945s$ (\textit{offline-online stabilization} without supremizer)\\
& $1730s$ (\textit{offline-only stabilization} with supremizer)\\
 \hline
\multirow{3}{*}{\textit{Online} time ($\mathbb{P}_1/{\mathbb{P}_1}$)} & $103s$ (\textit{offline-online stabilization} with supremizer)\\
& $82s$ (\textit{offline-online stabilization} without supremizer)\\
& $81s$ (\textit{offline-only stabilization} with supremizer)\\
 \hline
\multirow{3}{*}{\textit{Online} time ($\mathbb{P}_2/{\mathbb{P}_2}$)} & $242s$ (\textit{offline-online stabilization} with supremizer)\\
& $180s$ (\textit{offline-online stabilization} without supremizer)\\
& $90s$ (\textit{offline-only stabilization} with supremizer)\\
\hline
\hline
\end{tabular}\vspace{0.3cm}
\captionof{table}{Stokes problem: Computational details of unsteady Stokes problem \eqref{eq:test1_unsteady}.}
\label{tab:comput_unsteady}
\end{center}

\subsection{Numerical results for $\mathbb{P}_1/{\mathbb{P}_0}$}
\label{subsubsec: P1PO_unsteady}
In this subsection we show some results for the error comparison between the different stabilization options using the lowest order FE pair $\mathbb{P}_1/{\mathbb{P}_0}.$ The choice of stabilization term in equation  \eqref{eq:unsteady_discrete_stab} for lowest order element is as follows \cite{QV}:
\begin{equation}
s^{p,q}_{h}(q_{h};\boldsymbol\mu):=\delta\sum_{\sigma\in\Gamma_h}h_{\sigma}\int_{\sigma}\left[p_{h}\right]_{\sigma}\left[q_{h}\right]_{\sigma},
\label{eq:P1P0 STAB}
\end{equation}
where $\Gamma_h$ is the set of all edges $\sigma$ of the triangulation except for those belonging to the boundary $\partial\Omega$, $h_{\sigma}$ is the length of $\sigma$ and $\left[q_{h}\right]_{\sigma}$ denotes its jump across $\sigma$.
\par
The motivation in doing this case is to support the \textit{offline-online stabilization}, i.e, we want to show, by doing different numerical experiments that the \textit{offline-online stabilization} is the best way to stabilize whatever the stabilization we chose. For instance, in subsection \ref{subsec: PkPk} we chose the Franca-Hughes stabilization, which has different stabilization terms as compared to this subsection.\par
\begin{figure}[H]
  \centering
  \begin{subfigure}{0.48\textwidth}
    \centering\includegraphics[width=\textwidth]{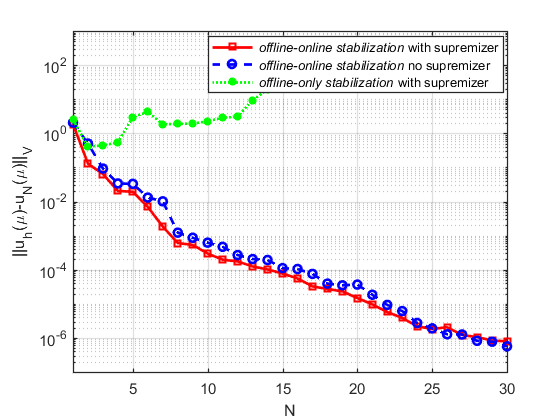}
  \end{subfigure}%
  \quad%
  \begin{subfigure}{0.48\textwidth}
    \centering\includegraphics[width=\textwidth]{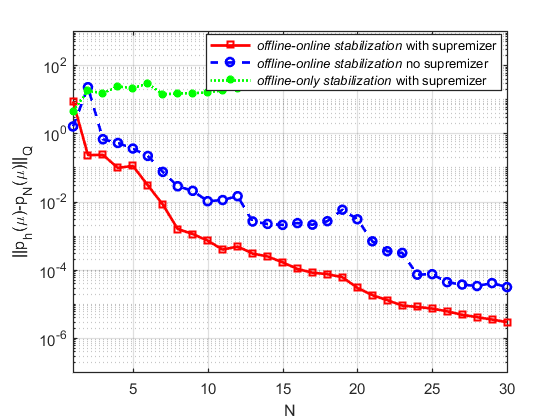}
  \end{subfigure}
  \caption{Stokes problem: $L^2$-error in time for velocity (left) and pressure (right) with stabilization coefficient $\delta=0.05$ and $\Delta{t}=0.02.$ using $\mathbb{P}_1/{\mathbb{P}_0}.$}
    \label{fig: L2U_P1P0}
\end{figure}

We plot the $L^2-$ error in time for velocity and pressure, respectively in Fig. \ref{fig: L2U_P1P0}. These results further strengthen our claim that the \textit{offline-online stabilization} is the best way to stabilize.

\subsection{Sensitivity on $\Delta{t}$}
\label{subsubsec: sensitivity_deltaT}
Consistently stabilized FE methods have complications while working with small time steps. These complications are reported in \cite{Max2004,Max2004a} and references therein. The analysis found in \cite{Max2004a} established that
\begin{center}
$\Delta t>\delta{h}^2$,
\end{center}
is a sufficient condition to avoid instabilities. Later on a detailed study and series of numerical experiments are performed in \cite{Max2007} and it is established that the fully discrete problem \eqref{eq:unsteady_discrete_stab} is conditionally stable with the condition
\begin{equation}
\Delta{t}/\delta{h}^2\geq\delta,
\label{eq: deltaT}
\end{equation}
where $\Delta{t}$ is the time step, $\delta$ is the stabilization coefficient independent of the spatial grid size $h.$\par
\begin{figure}[H]
  \centering
  \begin{subfigure}{0.48\textwidth}
    \centering\includegraphics[width=\textwidth]{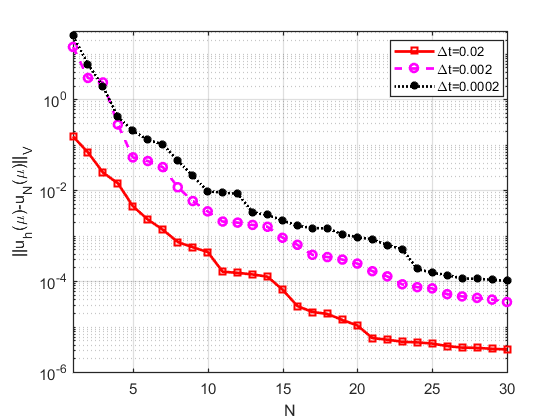}
  \end{subfigure}%
  \quad%
  \begin{subfigure}{0.48\textwidth}
    \centering\includegraphics[width=\textwidth]{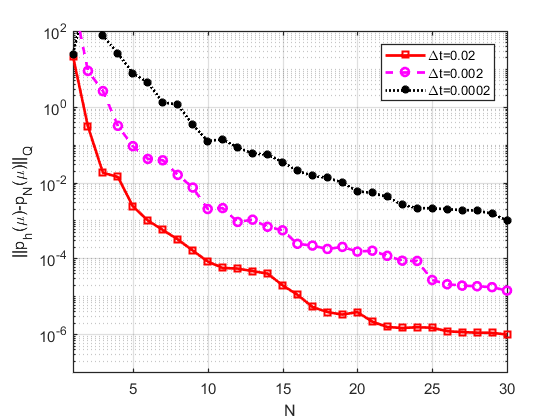}
  \end{subfigure}
  \caption{Stokes problem: Franca-Hughes stabilization; $L^2$-error in time for Velocity (left) and pressure (right) using $\mathbb{P}_2/\mathbb{P}_2$ and $\delta=0.05$, $\Delta{t}=0.02,0.002,0.0002$. \textit{offline-online stabilization} without supremizer.}
    \label{fig: U_sensitive}
\end{figure}

In this subsection we present some numerical results to see the variation of $\Delta{t}$ on the error between FE and RB solutions. We use the \textit{offline-online stabilization} without supremizer to plot the error between FE and RB solution for velocity (left) and pressure (right) in Fig. \ref{fig: U_sensitive}. We fix the value of stabilization coefficient $\delta=0.05$ \par
From these error plots, we observe that $\Delta t=0.02$ (in this case, not generally) is the best value. If we decrease the value of $\Delta t$, keeping $\delta$ and $h$ fixed, i.e, we are decreasing the left hand side of \eqref{eq: deltaT}, which increases the error.

\section{Parametrized unsteady Navier-Stokes problem}
\label{sec: Unsteady_NS}
In this section, we develop a stabilized RB method using SUPG stabilization method for the approximation of unsteady Navier-Stokes problem in reduced order parametric setting. Let $\Omega\subset\mathbb{R}^2$, be a reference configuration, and we assume that current configuration $\Omega_o(\boldsymbol\mu)$ can be obtained as the image of map $\boldsymbol{T}(.;\boldsymbol\mu):\mathbb{R}^2\rightarrow\mathbb{R}^2,$ i.e. $\Omega_o(\boldsymbol\mu)=\boldsymbol{T}(\Omega;\boldsymbol\mu).$
First we define the unsteady Navier-Stokes problem on a domain $\Omega_o(\boldsymbol\mu)$ in $\mathbb{R}^{2}$. We consider the fluid flow in a region $\Omega_o(\boldsymbol\mu)$, bounded by walls and driven by a body force $\boldsymbol{f}(\boldsymbol\mu)$. The fluid velocity and pressure are the functions $\boldsymbol{u_o}(t;\boldsymbol\mu)$ for $\boldsymbol\mu\in\mathbb{P}, 0\leq{t}\leq T$ and $p_o(t;\boldsymbol\mu)$ for $0<t\leq T$, respectively which satisfies
\begin{equation}
\begin{cases}
\dfrac{\partial\boldsymbol{u_o}}{\partial{t}}-\nu\Delta{\boldsymbol{u}_o}+\boldsymbol{u_o}\cdot\nabla{\boldsymbol{u_o}}+\nabla{p_o}= \boldsymbol{f}(\boldsymbol\mu) & \text{ in } \Omega_o(\boldsymbol\mu)\times\left(0,T\right),\  \\
\rm div {\boldsymbol{u}_o}=0 & \text{ in } \Omega_o(\boldsymbol\mu)\times\left(0,T\right),\\
\boldsymbol{u}_o=\boldsymbol{g} & \text{ on } \partial\Omega\times\left(0,T\right),\\
\boldsymbol{u}_o|_{t=0}=\boldsymbol{u_0} & \text{ in } \Omega_o(\boldsymbol\mu).
\end{cases}
\label{eq:unsteady_original_NS}
\end{equation}
By multiplying \eqref{eq:unsteady_original_NS} with velocity and pressure test functions $\boldsymbol{v}$ and $q$, respectively, integrating by parts, and tracing everything back onto the reference domain $\Omega,$ we obtain the following parametrized weak formulation of \eqref{eq:unsteady_original_NS}:\par
for a given $\boldsymbol\mu\in\mathbb{P},$ find $\boldsymbol{u}(t;\boldsymbol\mu)\in\boldsymbol{V}$ and $p(t;\boldsymbol\mu))\in Q$ such that
\begin{equation}
\begin{cases}
m(\dfrac{\partial\boldsymbol{u}}{\partial{t}},\boldsymbol{v};\boldsymbol\mu)+a(\boldsymbol{u},\boldsymbol{v};\boldsymbol\mu)+c(\boldsymbol{u},\boldsymbol{u},\boldsymbol{v};\boldsymbol\mu)+b(\boldsymbol{v},p;\boldsymbol\mu)=F(\boldsymbol{v};\boldsymbol\mu) & \forall \, v\in \boldsymbol{V},t>0,\\
b(\boldsymbol{u},q;\boldsymbol\mu)=G(q;\boldsymbol\mu) & \forall \, q\in Q,t>0,\\
\boldsymbol{u}|_{t=0}=\boldsymbol{u_0},
\end{cases}
\label{eq:unsteady_reference_NS}
\end{equation}
where the bilinear forms are given in \eqref{eq: unsteady_bilinear_forms} and trilinear form is defined as:
\begin{equation}
c(\boldsymbol{u},\boldsymbol{v},\boldsymbol{w};\boldsymbol\mu)=\int_{\Omega}u_i \chi_{ji}(x;\boldsymbol\mu)\dfrac{\partial{v_m}}{\partial x_j}w_m d\boldsymbol{x}.
\label{eq: unsteady_forms_NS}
\end{equation}
The tensors $\boldsymbol\kappa$, $\boldsymbol\chi$ and scalar $\pi$ are given by \eqref{eq: unsteady_tensors}.
\subsection{Discrete Finite Element formulation}
As in the previous part for unsteady Stokes problem, let us now discretize problem \eqref{eq:unsteady_reference_NS}. Consider $\lbrace T_h\rbrace_{h>0}$ be the triangulations and $h$ denotes a discretization parameter \cite{Raviart1986,Max1989}. Let $\boldsymbol{V}_h$ and $Q_h$ be two finite dimensional spaces such that $\boldsymbol{V}_h\subset\boldsymbol{H}^1(\Omega)$ and $Q_h\subset L^2_0(\Omega)$. We use implicit Euler scheme for time derivative term. We consider a partition of the interval $[0,T]$ into $K$ sub-intervals of equal length $\Delta{t}=T/K$ and $t^k=k\Delta{t}, 0\leq{k}\leq{K}.$ We approximate the time derivative in the $(k)-th$ time layer as
\begin{equation}
\dfrac{\partial\boldsymbol{u}_{h}(t^k)}{\partial{t}}\approx \dfrac{\boldsymbol{u}_{h}^{k}-\boldsymbol{u}_{h}^{k-1}}{\Delta{t}},
\label{eq:euler}
\end{equation}
where $\Delta t$ is a constant time step. We define the semi discrete FE approximation problem of \eqref{eq:unsteady_reference_NS} while  using \eqref{eq:euler} in \eqref{eq:unsteady_reference_NS} we get as follows:\\
for a given $\boldsymbol\mu\in\mathbb{P},$ and $(\boldsymbol{u}_h^{k-1}(\boldsymbol\mu),p_h^{k-1}(\boldsymbol\mu)),$ find $\boldsymbol{u}_h^{k}(t;\boldsymbol\mu)\in\boldsymbol{V}_h$ and $p_h^{k}(t;\boldsymbol\mu)\in Q_h$ such that
\begin{equation}
\begin{cases}
\frac{1}{\Delta{t}}m(\boldsymbol{u}_h^{k},\boldsymbol{v}_h;\boldsymbol\mu)+a(\boldsymbol{u}_h^{k},\boldsymbol{v}_h;\boldsymbol\mu)+c(\boldsymbol{u}_h^{k},\boldsymbol{u}_h^{k},\boldsymbol{v}_h;\boldsymbol\mu)\\+b(\boldsymbol{v}_h,p_h^{k};\boldsymbol\mu)=F(\boldsymbol{v}_h;\boldsymbol\mu)+\frac{1}{\Delta{t}}m(\boldsymbol{u}_h^{k-1},\boldsymbol{v}_h;\boldsymbol\mu) & \forall \, \boldsymbol{v}_h\in \boldsymbol{V}_h,\\
b(\boldsymbol{u}_h^{k},q_h;\boldsymbol\mu)=G(q_h;\boldsymbol\mu) & \forall \, q_h\in Q_h,\\
\boldsymbol{u}_{h}^{0}=\boldsymbol{u}_{0,h}.
\end{cases}
\label{eq:unsteady_discrete_time_NS}
\end{equation}
The algebraic formulation of \eqref{eq:unsteady_discrete_time_NS} can be written as:
\begin{equation}
\begin{split}
\left[\begin{array}{cc} \dfrac{M(\boldsymbol\mu)}{\Delta t}+ A(\boldsymbol\mu)+C(\boldsymbol{u}(t^{k};\boldsymbol\mu);\boldsymbol\mu)& B^T(\boldsymbol\mu)\\ B(\boldsymbol\mu) & \boldsymbol{0} \end{array}\right]
\left[\begin{array}{c} \boldsymbol{U}(t^{k};\boldsymbol\mu) \\ \boldsymbol{P}(t^{k};\boldsymbol\mu) \end{array}\right]
= \left[\begin{array}{c} \boldsymbol{\bar{f}}(\boldsymbol\mu) \\ \boldsymbol{\bar{g}}(\boldsymbol\mu) \end{array}\right] \\+
\left[\begin{array}{cc} \dfrac{M(\boldsymbol\mu)}{\Delta t} & \boldsymbol{0}\\ \boldsymbol{0} & \boldsymbol{0} \end{array}\right]
\left[\begin{array}{c} \boldsymbol{U}(t^{k-1};\boldsymbol\mu) \\ \boldsymbol{P}(t^{k-1};\boldsymbol\mu) \end{array}\right],
\end{split}
\label{eq: Algebraic_unsteady_time_NS}
\end{equation}
where the matrices corresponding to bilinear forms, and the vectors are given in \eqref{eq: matrices_vectors}. The matrix corresponding to nonlinear form is defined as:
\begin{equation}
\left(C(\boldsymbol{u}(t;\boldsymbol\mu);\boldsymbol\mu)\right)_{ij}=\sum_{m=1}^{\mathcal{N}_u}\boldsymbol{u}_h^{m}(t;\boldsymbol\mu)c(\boldsymbol\phi_m^h,\boldsymbol\phi_j^h,\boldsymbol\phi_i^h;\boldsymbol\mu),
\label{eq: matrices_vectors_NS}
\end{equation}
where $\boldsymbol\phi_i^h$ and $\boldsymbol\psi_j^h$, are the basis functions of $\boldsymbol{V}_h$ and $Q_h$ respectively. As in previous case, we impose the affine parametric dependence on these matrices and vectors and we skip the detail here.
\label{subsec:Unsteady_NS_Discrete}
\subsection{Stabilized Finite Element formulation}
In this section we give the stabilized formulation of time-dependent Navier-Stokes equations defined in previous section. We use the SUPG stabilization method \cite{Brooks1982} first in full order, and then, we project on reduced spaces to fulfill the reduced inf-sup condition.
\par
The stabilized FE formulation of \eqref{eq:unsteady_reference_NS} read as:
for a given $\boldsymbol\mu\in\mathbb{P},$ find $\boldsymbol{u}(t;\boldsymbol\mu)\in\boldsymbol{V}$ and $p(t;\boldsymbol\mu))\in Q$ such that
\begin{equation}
\begin{cases}
\frac{1}{\Delta{t}}m(\boldsymbol{u}_h^{k},\boldsymbol{v}_h;\boldsymbol\mu)+a(\boldsymbol{u}_h^{k},\boldsymbol{v}_h;\boldsymbol\mu)+c(\boldsymbol{u}_h^{k},\boldsymbol{u}_h^{k},\boldsymbol{v}_h;\boldsymbol\mu)\\+b(\boldsymbol{v}_h,p_h^{k};\boldsymbol\mu)-s^{ut,v}_{h}(\boldsymbol{u}_{h},\boldsymbol{v}_{h};\boldsymbol\mu)-s^{p,v}_{h}(p_h,\boldsymbol{v}_{h};\boldsymbol\mu)=F(\boldsymbol{v}_h;\boldsymbol\mu)+\frac{1}{\Delta{t}}m(\boldsymbol{u}_h^{k-1},\boldsymbol{v}_h;\boldsymbol\mu) & \forall \, \boldsymbol{v}_h\in \boldsymbol{V}_h,\\
b(\boldsymbol{u}_h^{k},q_h;\boldsymbol\mu)-s^{ut,q}_{h}(\boldsymbol{u}_{h},q_{h};\boldsymbol\mu)-s^{p,q}_{h}(p_h,q_{h};\boldsymbol\mu)=G(q_h;\boldsymbol\mu) & \forall \, q_h\in Q_h,\\
\boldsymbol{u}_{h}^{0}=\boldsymbol{u}_{0,h}.
\end{cases}
\label{eq:unsteady_reference_NS_stab}
\end{equation}
where $s^{ut,v}_{h}(.,.;\boldsymbol\mu)$, $s^{p,v}_{h}(.,.;\boldsymbol\mu)$, $s^{ut,q}_{h}(.,.;\boldsymbol\mu)$ and $s^{p,q}_{h}(.,.;\boldsymbol\mu)$ are the stabilization terms \cite{QV} defined as:
\begin{equation}
\begin{split}
s^{ut,v}_{h}(\boldsymbol{u}_h,\boldsymbol{v}_h;\boldsymbol\mu)&:=\delta\sum_{K}h_K^{2}\int_K(\dfrac{\partial\boldsymbol{u}}{\partial{t}}-\nu\Delta\boldsymbol{u}_{h}+\boldsymbol{u}_h\cdot\nabla{\boldsymbol{u}_h},\boldsymbol{u}_h\cdot\nabla{\boldsymbol{v}_h}),
\end{split}
\label{eq: stab_terms_un_NS}
\end{equation}

\begin{equation}
\begin{split}
s^{p,v}_{h}(p_h,\boldsymbol{v}_h;\boldsymbol\mu)&:=\delta\sum_{K}h_K^{2}\int_K(\nabla{p}_{h},\boldsymbol{u}_h\cdot\nabla{\boldsymbol{v}_h}),
\end{split}
\label{eq: stab_terms_un_NS_1}
\end{equation}

\begin{equation}
\begin{split}
s^{ut,q}_{h}(\boldsymbol{u}_h,q_h;\boldsymbol\mu)&:=\delta\sum_{K}h_K^{2}\int_K(\dfrac{\partial\boldsymbol{u}}{\partial{t}}-\nu\Delta\boldsymbol{u}_{h}+\boldsymbol{u}_h\cdot\nabla{\boldsymbol{u}_h},\nabla{q}_{h}),
\end{split}
\label{eq: stab_terms_un_NS_2}
\end{equation}

\begin{equation}
s^{p,q}_{h}(p_h,q_h;\boldsymbol\mu):=\delta\sum_{K}h_K^{2}\int_K(\nabla{p}_{h},\nabla{q}_{h}),
\label{eq: stab_terms_un_NS_3}
\end{equation}
where $\delta$ is the stabilization coefficient.
The stabilized algebraic formulation of \eqref{eq:unsteady_reference_NS_stab} reads as:
\begin{equation}
\begin{split}
\left[\begin{array}{cc} \dfrac{M(\boldsymbol\mu)}{\Delta t}+ A(\boldsymbol\mu)+\tilde{C}(\boldsymbol{u}(t^{k};\boldsymbol\mu);\boldsymbol\mu)& \tilde{B}^T(\boldsymbol\mu)\\ \tilde{B}(\boldsymbol\mu)+\dfrac{\tilde{M}(\boldsymbol\mu)}{\Delta t} & -S(\boldsymbol\mu) \end{array}\right]
\left[\begin{array}{c} \boldsymbol{U}(t^{k};\boldsymbol\mu) \\ \boldsymbol{P}(t^{k};\boldsymbol\mu) \end{array}\right]
= \left[\begin{array}{c} \boldsymbol{\bar{f}}(\boldsymbol\mu) \\ \boldsymbol{\tilde{\bar{g}}}(\boldsymbol\mu) \end{array}\right] \\+
\left[\begin{array}{cc} \dfrac{M(\boldsymbol\mu)}{\Delta t} & \boldsymbol{0}\\ \dfrac{\tilde{M}(\boldsymbol\mu)}{\Delta t} & \boldsymbol{0} \end{array}\right]
\left[\begin{array}{c} \boldsymbol{U}(t^{k-1};\boldsymbol\mu) \\ \boldsymbol{P}(t^{k-1};\boldsymbol\mu) \end{array}\right],
\end{split}
\label{eq: Algebraic_unsteady_time_NS_stab}
\end{equation}
where $\tilde{B}, \tilde{B}^T, \tilde{M}$ and $\tilde{C}$, are the sum of original matrices in formulation \eqref{eq: Algebraic_unsteady_time_NS} and the SUPG stabilization matrices. Similarly $\boldsymbol{\tilde{\bar{f}}}$ and $\boldsymbol{\tilde{\bar{g}}}$ are vectors on right hand side which are sum of original vectors in formulation \eqref{eq: Algebraic_unsteady_time_NS} and SUPG stabilization terms \cite{Formaggia2012}. These matrices and vectors can be written similar to Stokes case \eqref{eq: Stokes_stab_matrices}.
\label{subsec:UN_NS_StabFE}

\subsection{Reduced Basis formulation}
A reduced order approximation of velocity and pressure field is obtained by means of Galerkin projection on the RB spaces $\boldsymbol{V}_N, Q_N$ and $\tilde{\boldsymbol{V}}_N$, defined in \eqref{eq:velocity}, \eqref{eq:pressure_t} and \eqref{eq: sup_velocity}, respectively.\par
In the \textit{online} stage, the resulting reduced order approximation of \eqref{eq:unsteady_discrete_time_NS} is as follows:
for any parameter $\boldsymbol\mu\in \mathbb{P},$ find $\boldsymbol{u}_N^{k}(t;\boldsymbol\mu)\in\boldsymbol{V}_N$ and $p_N^{k}(t;\boldsymbol\mu)\in Q_N$ such that
\begin{equation}
\begin{cases}
\frac{1}{\Delta{t}}m(\boldsymbol{u}_N^{k},\boldsymbol{v}_N;\boldsymbol\mu)+a(\boldsymbol{u}_N^{k},\boldsymbol{v}_N;\boldsymbol\mu)+c(\boldsymbol{u}_N^{k},\boldsymbol{u}_N^{k},\boldsymbol{v}_N;\boldsymbol\mu)\\+b(\boldsymbol{v}_N,p_N^{k};\boldsymbol\mu)=F(\boldsymbol{v}_N;\boldsymbol\mu)+\frac{1}{\Delta{t}}m(\boldsymbol{u}_N^{k-1},\boldsymbol{v}_N;\boldsymbol\mu) & \forall \, \boldsymbol{v}_N\in \boldsymbol{V}_N,\\
b(\boldsymbol{u}_N^{k},q_N;\boldsymbol\mu)=G(q_N;\boldsymbol\mu) & \forall \, q_N\in Q_N,\\
\boldsymbol{u}_{N}^{0}=\boldsymbol{u}_{0,N}.
\end{cases}
\label{eq:unsteady_discrete_time_NS_RB}
\end{equation}
The algebraic formulation of \eqref{eq:unsteady_discrete_time_NS_RB} can be written as
\begin{equation}
\begin{split}
\left[\begin{array}{cc} \dfrac{M_N(\boldsymbol\mu)}{\Delta t}+ A_N(\boldsymbol\mu)+C_N(\boldsymbol{u}_N(t^{k};\boldsymbol\mu);\boldsymbol\mu)& B_N^T(\boldsymbol\mu)\\ B_N(\boldsymbol\mu) & \boldsymbol{0} \end{array}\right]
\left[\begin{array}{c} \boldsymbol{U}_N(t^{k};\boldsymbol\mu) \\ \boldsymbol{P}_N(t^{k};\boldsymbol\mu) \end{array}\right]
= \left[\begin{array}{c} \boldsymbol{\bar{f}}_N(\boldsymbol\mu) \\ \boldsymbol{\bar{g}}_N(\boldsymbol\mu) \end{array}\right] \\+\left[\begin{array}{cc} \dfrac{M_N(\boldsymbol\mu)}{\Delta t} & \boldsymbol{0}\\ \boldsymbol{0} & \boldsymbol{0} \end{array}\right]
\left[\begin{array}{c} \boldsymbol{U}_N(t^{k-1};\boldsymbol\mu) \\ \boldsymbol{P}_N(t^{k-1};\boldsymbol\mu) \end{array}\right],
\end{split}
\label{eq: Algebraic_unsteady_time_NS_RB}
\end{equation}
where, the reduced order matrices are defined as:
\begin{equation}
\begin{split}
M_N(\boldsymbol\mu)=Z^T_{u,s}M(\boldsymbol\mu)_{u,s}, \quad
A_N(\boldsymbol\mu)=Z^T_{u,s}A(\boldsymbol\mu)Z_{u,s}, \quad B_N(\boldsymbol\mu)=Z^T_{p}B(\boldsymbol\mu)Z_{u}, \\
C_N(.;\boldsymbol\mu)=Z^T_{u}C(.;\boldsymbol\mu)Z_{u}, \quad  \boldsymbol{\bar{f}}_N(\boldsymbol\mu)=Z^T_{u}\boldsymbol{\bar{f}}(\boldsymbol\mu),\quad
\boldsymbol{\bar{g}}_N(\boldsymbol\mu)=Z^T_{p}\boldsymbol{\bar{g}}(\boldsymbol\mu).
\end{split}
\label{eq:UN_NS_RBmatrices}
\end{equation}
\label{subsec:UN_NS_RB}
\subsection{Stabilized Reduced Basis formulation}
We write the stabilized formulation of \eqref{eq:unsteady_discrete_time_NS_RB} as follows:
for any parameter $\boldsymbol\mu\in \mathbb{P},$ find $\boldsymbol{u}_N^{k}(t;\boldsymbol\mu)\in\boldsymbol{V}_N$ and $p_N^{k}(t;\boldsymbol\mu)\in Q_N$ such that
\begin{equation}
\begin{cases}
\frac{1}{\Delta{t}}m(\boldsymbol{u}_N^{k},\boldsymbol{v}_N;\boldsymbol\mu)+a(\boldsymbol{u}_N^{k},\boldsymbol{v}_N;\boldsymbol\mu)+c(\boldsymbol{u}_N^{k},\boldsymbol{u}_N^{k},\boldsymbol{v}_N;\boldsymbol\mu)\\+b(\boldsymbol{v}_N,p_N^{k};\boldsymbol\mu)-s^{ut,v}_{N}(\boldsymbol{u}_{N},\boldsymbol{v}_{N};\boldsymbol\mu)-s^{p,v}_{N}(p_N,\boldsymbol{v}_{N};\boldsymbol\mu)=F(\boldsymbol{v}_N;\boldsymbol\mu)+\frac{1}{\Delta{t}}m(\boldsymbol{u}_N^{k-1},\boldsymbol{v}_N;\boldsymbol\mu) & \forall \, \boldsymbol{v}_N\in \boldsymbol{V}_N,\\
b(\boldsymbol{u}_N^{k},q_N;\boldsymbol\mu)-s^{ut,q}_{N}(\boldsymbol{u}_{N},q_{N};\boldsymbol\mu)-s^{p,q}_{N}(p_N,q_{N};\boldsymbol\mu)=G(q_N;\boldsymbol\mu) & \forall \, q_N\in Q_N,\\
\boldsymbol{u}_{N}^{0}=\boldsymbol{u}_{0,N}.
\end{cases}
\label{eq:unsteady_discrete_time_NS_RB_stab}
\end{equation}
where $s^{ut,v}_{N}(.,.;\boldsymbol\mu)$, $s^{p,v}_{N}(.,.;\boldsymbol\mu)$, $s^{ut,q}_{N}(.,.;\boldsymbol\mu)$ and $s^{p,q}_{N}(.,.;\boldsymbol\mu)$ are reduced order stabilization terms defined as:
\begin{equation}
s^{ut,v}_{N}(\boldsymbol{u}_N,\boldsymbol{v}_N;\boldsymbol\mu):=\delta\sum_{K}h_K^{2}\int_K(\dfrac{\partial\boldsymbol{u}}{\partial{t}}-\nu\Delta\boldsymbol{u}_{N}+\boldsymbol{u}_N\cdot\nabla{\boldsymbol{u}_N},\boldsymbol{u}_N\cdot\nabla{\boldsymbol{v}_N}),
\label{eq: stab_terms_un_NS_RB}
\end{equation}

\begin{equation}
s^{p,v}_{N}(p_N,\boldsymbol{v}_N;\boldsymbol\mu):=\delta\sum_{K}h_K^{2}\int_K(\nabla{p}_{N},\boldsymbol{u}_N\cdot\nabla{\boldsymbol{v}_N}),
\label{eq: stab_terms_un_NS_1_RB}
\end{equation}

\begin{equation}
s^{ut,q}_{N}(\boldsymbol{u}_N,q_N;\boldsymbol\mu):=\delta\sum_{K}h_K^{2}\int_K(\dfrac{\partial\boldsymbol{u}}{\partial{t}}-\nu\Delta\boldsymbol{u}_{N}+\boldsymbol{u}_N\cdot\nabla{\boldsymbol{u}_N},\nabla{q}_{N}),
\label{eq: stab_terms_un_NS_2_RB}
\end{equation}

\begin{equation}
s^{p,q}_{h}(p_N,q_N;\boldsymbol\mu):=\delta\sum_{K}h_K^{2}\int_K(\nabla{p}_{N},\nabla{q}_{N}),
\label{eq: stab_terms_un_NS_3_RB}
\end{equation}
The algebraic formulation of \eqref{eq:unsteady_discrete_time_NS_RB_stab} can be written as:
\begin{equation}
\begin{split}
\left[\begin{array}{cc} \dfrac{M_N(\boldsymbol\mu)}{\Delta t}+ A_N(\boldsymbol\mu)+\tilde{C}_N(\boldsymbol{u}_N(t^{k};\boldsymbol\mu);\boldsymbol\mu)& \tilde{B}_N^T(\boldsymbol\mu)\\ \tilde{B}_N(\boldsymbol\mu)+\dfrac{\tilde{M}_N(\boldsymbol\mu)}{\Delta t} & -S_N(\boldsymbol\mu) \end{array}\right]
\left[\begin{array}{c} \boldsymbol{U}_N(t^{k};\boldsymbol\mu) \\ \boldsymbol{P}_N(t^{k};\boldsymbol\mu) \end{array}\right]
= \left[\begin{array}{c} \boldsymbol{\bar{f}}_N(\boldsymbol\mu) \\ \boldsymbol{\tilde{\bar{g}}}_N(\boldsymbol\mu) \end{array}\right] \\+
\left[\begin{array}{cc} \dfrac{M_N(\boldsymbol\mu)}{\Delta t} & \boldsymbol{0}\\ \dfrac{\tilde{M}_N(\boldsymbol\mu)}{\Delta t} & \boldsymbol{0} \end{array}\right]
\left[\begin{array}{c} \boldsymbol{U}_N(t^{k-1};\boldsymbol\mu) \\ \boldsymbol{P}_N(t^{k-1};\boldsymbol\mu) \end{array}\right],
\end{split}
\label{eq: Algebraic_unsteady_time_NS_stabRB}
\end{equation}
where $\tilde{B}_N, \tilde{B}_N^T, \tilde{M}_N$ and $\tilde{C}_N$ are RB stabilized matrices, and can be obtained similarly as \eqref{eq:UN_NS_RBmatrices}.\par
\label{section:NS_UNSTEADY}

\section{Numerical results and discussion}
\label{sec:Results_UNS}
In this section we apply the stabilized RB model for unsteady Navier-Stokes problem presented in section \ref{sec: Unsteady_NS} and subsections therein to \textit{lid-driven cavity} flow problem on parametrized domain shown in Fig. \ref{fig:Domain}.
We first show some numerical results for only physical parameterization in subsection \ref{subsec:physical_results_UNS}, and then, we show numerical results for both physical and geometrical parametrization in subsection \ref{subsec:geo_results_UNS}. In both cases we compare and discuss the three options; $(i)$ \textit{offline-online stabilization} with supremizer, $(ii)$ \textit{offline-online stabilization} without supremizer, $(iii)$ \textit{offline-only stabilization} with supremizer.
\label{sec:Results_UNS}
\subsection{Results for physical parameter case only}
\label{subsec:physical_results_UNS}
The parameter in this case is only the physical parameter, i.e, the Reynolds number and is denoted by $\mu$. The details of computation is summarized in Table \ref{tab:comput_steady_UN_NS}.

\begin{figure}[H]
  \centering
  \begin{subfigure}{0.48\textwidth}
    \centering\includegraphics[width=\textwidth]{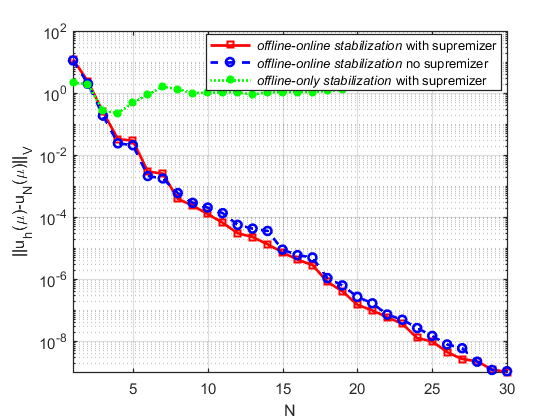}
  \end{subfigure}%
  \quad%
  \begin{subfigure}{0.48\textwidth}
    \centering\includegraphics[width=\textwidth]{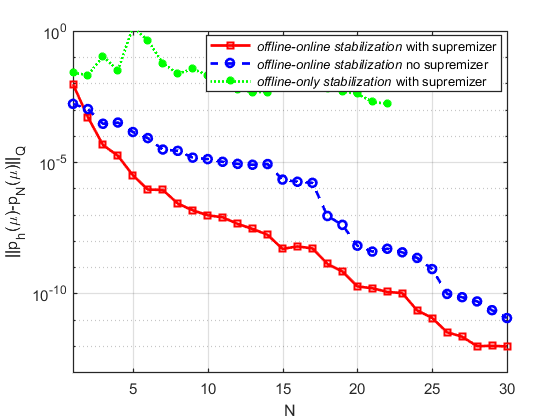}
  \end{subfigure}
  \caption{Navier-Stokes problem with SUPG stabilization; physical parametrization on cavity flow; Error between FE and RB solution for velocity (left) and pressure (right) using $\mathbb{P}_1/{\mathbb{P}_1}$.}
    \label{fig:UP_SUPG_tP1}
\end{figure}
\begin{figure}[H]
  \centering
  \begin{subfigure}{0.48\textwidth}
    \centering\includegraphics[width=\textwidth]{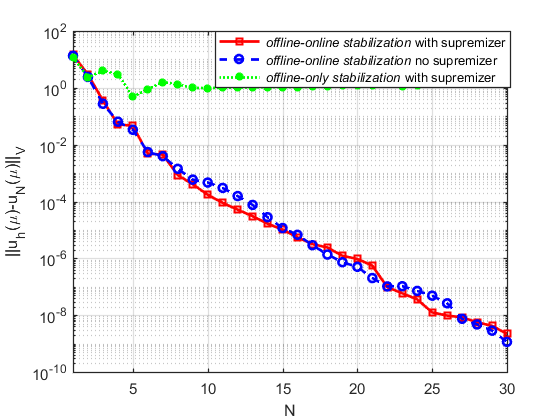}
  \end{subfigure}%
  \quad%
  \begin{subfigure}{0.48\textwidth}
    \centering\includegraphics[width=\textwidth]{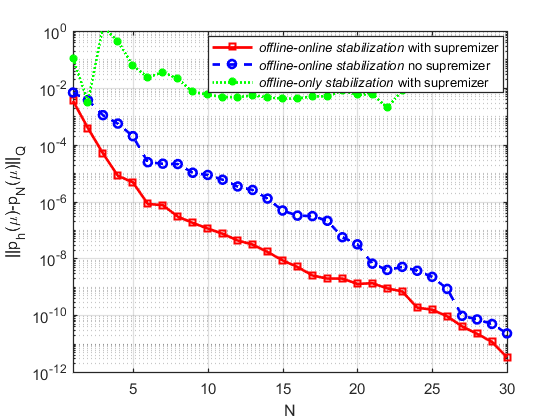}
  \end{subfigure}
  \caption{Navier-Stokes problem with SUPG stabilization; physical parametrization on cavity flow; Error between FE and RB solution for velocity (left) and pressure (right) using $\mathbb{P}_2/{\mathbb{P}_2}$.}
    \label{fig:UP_SUPG_tP2}
\end{figure}

\begin{center}
 \begin{tabular}{|| l | l||}
 \hline\hline
Physical parameter & $\mu$ (Reynolds number)\\
 \hline
Range of $\mu$  & [100,200]  \\
 \hline
\textit{Online} $\mu$ (example)  & 130  \\
 \hline
FE degrees of freedom & 5934 ($\mathbb{P}_1/{\mathbb{P}_1}$)\\
 \hline
RB dimension & $N_u=N_s=N_p=30$\\
 \hline
\multirow{2}{*}{\textit{Offline} time ($\mathbb{P}_1/{\mathbb{P}_1}$)} & $40612s$ (\textit{offline-online stabilization} with supremizer)\\
& $38781s$ (\textit{offline-online stabilization} without supremizer)\\
 \hline
\multirow{2}{*}{\textit{Online} time ($\mathbb{P}_1/{\mathbb{P}_1}$)} & $4640s$ (\textit{offline-online stabilization} with supremizer)\\
& $4040s$ (\textit{offline-online stabilization} without supremizer)\\
 \hline
Time step & $0.02$\\
 \hline
Final time & $0.5$\\
\hline\hline
\end{tabular}\vspace{0.3cm}
\captionof{table}{Navier-Stokes problem with physical parameter only: Computational details of unsteady Navier-Stokes problem without Empirical Interpolation.}
\label{tab:comput_steady_UN_NS}
\end{center}
Figure \ref{fig:UP_SUPG_tP1} plots the $L^2$-error in time for velocity (left) and pressure (right) using $\mathbb{P}_1/{\mathbb{P}_1}$ FE pair. Similarly results for velocity and pressure using $\mathbb{P}_2/{\mathbb{P}_2}$ FE pair are shown in Fig. \ref{fig:UP_SUPG_tP2}. In all numerical results presented in this section, we observe that the \textit{offline-online stabilization} without supremizer has better performance for velocity in terms of error. However, in case of pressure, our results show that supremizer is still improving the error but on the other hand addition of supremizer is computationally expensive. The \textit{offline-only stabilization} is not accurate also in this case.

\subsection{Results for physical and geometrical parameters}
\label{subsec:geo_results_UNS}
In this section we present some numerical results for unsteady Navier Stokes problem with physical and geometrical parameters. The computation details are presented in Table \ref{tab:comput_steady_NS_unsteady}. We recall that we are not using any ``hyper-reduction" technique to improve online performance at the moment. Our interest at the moment is in a preliminary testing of accuracy and stability.
\par
\begin{figure}[H]
  \centering
  \begin{subfigure}{0.48\textwidth}
    \centering\includegraphics[width=\textwidth]{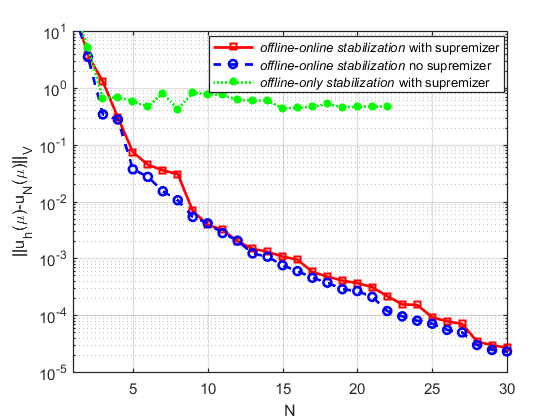}
  \end{subfigure}%
  \quad%
  \begin{subfigure}{0.48\textwidth}
    \centering\includegraphics[width=\textwidth]{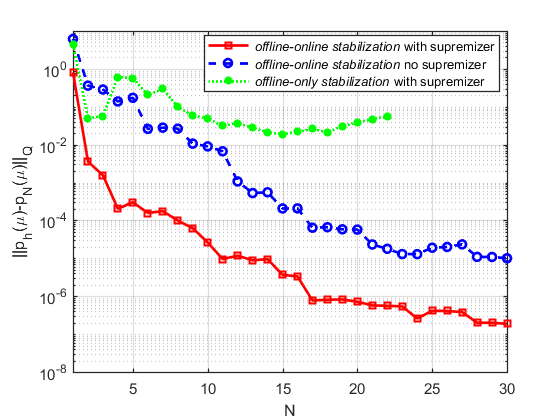}
  \end{subfigure}
  \caption{Navier-Stokes problem with SUPG stabilization using $\mathbb{P}_1/{\mathbb{P}_1}$: Velocity (left) and pressure (right) error for physical and geometrical parameters on cavity flow.}
    \label{fig:L2errorU_P2P2_NSGEO}
\end{figure}

Figure \ref{fig:L2errorU_P2P2_NSGEO} illustrates the error between FE and RB solution for velocity (left) and pressure (right) using $\mathbb{P}_1/{\mathbb{P}_1}$ FE pair. We observe that the error between two solutions, obtained by using \textit{offline-online stabilization} with/without supremizer is negligible in case of velocity. However, in case of pressure, supremizer has better performance. We have similar results for  $\mathbb{P}_2/{\mathbb{P}_2}$ FE pair that we do not show here.


\begin{center}
 \begin{tabular}{|| l | l||}
 \hline\hline
Physical parameter & $\mu_1$ (Reynolds number)\\
 \hline
Geometrical parameter & $\mu_2$ (horizontal length of domain)\\
 \hline
Range of $\mu_1$  & [100,200]  \\
 \hline
Range of $\mu_2$  & [1.5,3]  \\
 \hline
$\mu_1$ \textit{online} (example)  & 130  \\
 \hline
$\mu_2$ \textit{online} (example)  & 2  \\
 \hline
FE degrees of freedom & 6222 ($\mathbb{P}_1/{\mathbb{P}_1}$)\\
 \hline
RB dimension & $N_u=N_s=N_p=30$\\
 \hline
\multirow{2}{*}{\textit{Offline} time ($\mathbb{P}_1/{\mathbb{P}_1}$)} & $44693s$ (\textit{offline-online stabilization} with supremizer)\\
& $40153s$ (\textit{offline-online stabilization} without supremizer)\\
 \hline
\multirow{2}{*}{\textit{Online} time ($\mathbb{P}_1/{\mathbb{P}_1}$)} & $5169s$ (\textit{offline-online stabilization} with supremizer)\\
& $4724s$ (\textit{offline-online stabilization} without supremizer)\\
 \hline
Time step & $0.02$\\
 \hline
Final time & $0.5$\\
\hline\hline
\end{tabular}\vspace{0.3cm}
\captionof{table}{Computational details for unsteady Navier-Stokes problem with physical and geometrical parameters: stabilization and computational reduction.}
\label{tab:comput_steady_NS_unsteady}
\end{center}
\nocite{siddiqui2015}

\section{Concluding remarks}
\label{sec:Conclusion_NStokes}
In this work we have developed a stabilized RB method for the approximation of unsteady parametrized Stokes and Navier-Stokes problem. We have extended the analysis carried out in our previous work \cite{Ali2018} to the unsteady problems. The RB formulation is built, using the classical residual based stabilization technique in full order during the \textit{offline} stage and, then, projecting on the RB space. We have compared our approach with the existing approaches based on supremizers \cite{Veroy2007} through numerical experiments. In particular, the comparison between \textit{offline-online stabilization} with/without supremizer and \textit{offline-only stabilization} for unsteady Stokes and Navier-Stokes problems is presented. Our results in this work are consistent with those of the steady Stokes and Navier-Stokes case \cite{Ali2018}. On the basis of numerical results the main observations are as it follows:
\begin{itemize}
\item \textit{offline-online stabilization} is the most appropriate way to perform RB stabilization (if needed) for unsteady Stokes and Navier-Stokes problems;
\item using residual based stabilization, velocity is still better using \textit{offline-online stabilization} (without supremizer) even if pressure is improved in its accuracy by the supremizer enrichment;
\item \textit{offline-only stabilization} is not accurate. As in \cite{Ali2018}, this is due to the lack of consistency between the full and reduced order schemes, which occurs when solving the stabilized system during the \textit{offline} stage and non-stabilized system during the \textit{online} stage;
\item in terms of CPU time, the Taylor-Hood FE pair ($\mathbb{P}_2/{\mathbb{P}_1}$) is more expensive than ($\mathbb{P}_1/{\mathbb{P}_1}$) stabilized but less expensive than ($\mathbb{P}_2/{\mathbb{P}_2}$) stabilized (see, for instance Table \ref{tab:comput_unsteady});
\end{itemize}
\section*{Acknowledgements}
This work has been supported by the European Union Funding for Research and Innovation -- Horizon 2020 Program -- in the framework of European Research Council Executive Agency: H2020 ERC CoG 2015 AROMA-CFD project 681447 ``Advanced Reduced Order Methods with Applications in Computational Fluid Dynamics''. We also acknowledge the INDAM-GNCS project ``Advanced intrusive and non-intrusive model order reduction techniques and applications''.

\bibliographystyle{abbrv}
\bibliography{article}

\begin{thebibliography}{10}

\bibitem{Ali2018}
S.~Ali, F.~Ballarin, and G.~Rozza.
\newblock Stabilized reduced basis methods for parametrized steady stokes and
  navier–stokes equations.
\newblock {\em Computers \& Mathematics with Applications}, 80(11):2399 --
  2416, 2020.

\bibitem{BallarinChaconDelgadoGomezRozza2020}
F.~Ballarin, T.~Chacón~Rebollo, E.~Delgado~Ávila, M.~Gómez~Mármol, and
  G.~Rozza.
\newblock Certified reduced basis {VMS}-smagorinsky model for natural
  convection flow in a cavity with variable height.
\newblock {\em Computers \& Mathematics with Applications}, 80(5):973--989,
  2020.

\bibitem{Ballarin2015}
F.~Ballarin, A.~Manzoni, A.~Quarteroni, and G.~Rozza.
\newblock Supremizer stabilization of {POD}-{G}alerkin approximation of
  parametrized steady incompressible {N}avier-{S}tokes equations.
\newblock {\em International Journal for Numerical Methods in Engineering},
  102(5):1136--1161, 2015.

\bibitem{Becker2001}
R.~Becker and M.~Braack.
\newblock A finite element pressure gradient stabilization for the stokes
  equations based on local projections.
\newblock {\em Calcolo}, 38(4):173--199, 2001.

\bibitem{Max2004}
P.~B. Bochev, M.~D. Gunzburger, , and R.~Lehoucq.
\newblock On stabilized finite element methods for transient problems with
  varying time scales.
\newblock {\em Proceedings of ECOMAS 2004}.

\bibitem{Max2007}
P.~B. Bochev, M.~D. Gunzburger, and R.~B. Lehoucq.
\newblock On stabilized finite element methods for the {S}tokes problem in the
  small time step limit.
\newblock {\em International Journal for Numerical Methods in Fluids},
  53(4):573--597, 2007.

\bibitem{Max2004a}
P.~B. Bochev, M.~D. Gunzburger, and J.~N. Shadid.
\newblock On inf-sup stabilized finite element methods for transient problems.
\newblock {\em Computer Methods in Applied Mechanics and Engineering},
  193(15):1471 -- 1489, 2004.

\bibitem{boffi2013mixed}
D.~Boffi, F.~Brezzi, and M.~Fortin.
\newblock {\em Mixed finite element methods and applications}, volume~44.
\newblock Springer, 2013.

\bibitem{Hughes1980}
A.~Brooks and T.~Hughes.
\newblock Streamline {U}pwind/{P}etrov-{G}alerkin methods for advection
  dominated flows.
\newblock {\em Third International Conference on Finite Element Methods in
  Fluid Flow}, 2 , Calgary, Canada, Calgary Univ., 1980.

\bibitem{Brooks1982}
A.~N. Brooks and T.~J. Hughes.
\newblock Streamline {U}pwind/{P}etrov-{G}alerkin formulations for convection
  dominated flows with particular emphasis on the incompressible
  {N}avier-{S}tokes equations.
\newblock {\em Computer Methods in Applied Mechanics and Engineering}, 32:199
  -- 259, 1982.

\bibitem{Burman2009}
E.~Burman and M.~Fern\'{a}ndez.
\newblock Galerkin finite element methods with symmetric pressure stabilization
  for the transient stokes equations: Stability and convergence analysis.
\newblock {\em SIAM Journal on Numerical Analysis}, 47(1):409--439, 2009.

\bibitem{Christie1976}
I.~Christie, D.~F. Griffiths, A.~R. Mitchell, and O.~C. Zienkiewicz.
\newblock Finite element methods for second order differential equations with
  significant first derivatives.
\newblock {\em International Journal for Numerical Methods in Engineering},
  10(6):1389--1396, 1976.

\bibitem{Simone2008}
S.~Deparis.
\newblock Reduced basis error bound computation of parameter-dependent
  {N}avier-{S}tokes equations by the natural norm approach.
\newblock {\em SIAM Journal on Numerical Analysis}, 46(4):2039--2067, 2008.

\bibitem{Simone2009}
S.~Deparis and G.~Rozza.
\newblock Reduced basis method for multi-parameter-dependent steady
  {N}avier-{S}tokes equations: Applications to natural convection in a cavity.
\newblock {\em Journal of Computational Physics}, 228(12):4359 -- 4378, 2009.

\bibitem{Douglas1989}
J.~J. Douglas and J.~Wang.
\newblock An absolutely stabilized finite element formulation for the {S}tokes
  problem.
\newblock {\em Mathematics of Computations}, 52(186):495 -- 508, 1989.

\bibitem{Formaggia2012}
L.~Formaggia, F.~Saleri, and A.~Veneziani.
\newblock {\em Solving Numerical PDEs: Problems, Applications, Exercises}.
\newblock Springer-Verlag Mailand, 2012.

\bibitem{FRANCA1992}
L.~P. Franca and S.~L. Frey.
\newblock Stabilized finite element methods: {II}. {T}he incompressible
  {N}avier-{S}tokes equations.
\newblock {\em Computer Methods in Applied Mechanics and Engineering},
  99(2):209 -- 233, 1992.

\bibitem{FRANCA1992(1)}
L.~P. Franca, S.~L. Frey, and T.~J. Hughes.
\newblock Stabilized finite element methods: {I}. application to the
  advective-diffusive model.
\newblock {\em Computer Methods in Applied Mechanics and Engineering},
  95(2):253 -- 276, 1992.

\bibitem{Raviart1986}
V.~Girault and P.-A. Raviart.
\newblock {\em Finite Element Methods for Navier-Stokes Equations}, volume~5.
\newblock Springer, 1986.

\bibitem{Max1989}
M.~Gunzburger.
\newblock {\em Finite Element Methods for Viscous Incompressible Flows},
  volume~5.
\newblock Academicr, 1989.

\bibitem{Haasdonk2013}
B.~Haasdonk.
\newblock Convergence rates of the {POD}-greedy method.
\newblock {\em ESAIM: M2AN}, 47(3):859--873, 2013.

\bibitem{HANSBO1990}
P.~Hansbo and A.~Szepessy.
\newblock A velocity-pressure streamline diffusion finite element method for
  the incompressible {N}avier-{S}tokes equations.
\newblock {\em Computer Methods in Applied Mechanics and Engineering},
  84(2):175 -- 192, 1990.

\bibitem{Heinrich1977}
J.~C. Heinrich, P.~S. Huyakorn, O.~C. Zienkiewicz, and A.~R. Mitchell.
\newblock An ‘upwind’ finite element scheme for two-dimensional convective
  transport equation.
\newblock {\em International Journal for Numerical Methods in Engineering},
  11(1):131--143, 1977.

\bibitem{Heinrich1977(1)}
J.~C. Heinrich and O.~C. Zienkiewicz.
\newblock Quadratic finite element schemes for two-dimensional
  convective-transport problems.
\newblock {\em International Journal for Numerical Methods in Engineering},
  11(12):1831--1844, 1977.

\bibitem{Ali2018b}
S.~Hijazi, S.~Ali, G.~Stabile, F.~Ballarin, and G.~Rozza.
\newblock {\em The Effort of Increasing Reynolds Number in Projection-Based
  Reduced Order Methods: From Laminar to Turbulent Flows}, pages 245--264.
\newblock Springer International Publishing, 2020.

\bibitem{Hughes1979}
T.~Hughes and A.~Brooks.
\newblock A multi-dimensioal upwind scheme with no crosswind diffusion.
\newblock {\em Finite Element Methods for Convection Dominated Flows, New York,
  U.S.A.,}, 34:19--35, 1979.

\bibitem{HUGHES1987}
T.~J. Hughes and L.~P. Franca.
\newblock A new finite element formulation for computational fluid dynamics:
  {VII}. {T}he {S}tokes problem with various well-posed boundary conditions:
  Symmetric formulations that converge for all velocity/pressure spaces.
\newblock {\em Computer Methods in Applied Mechanics and Engineering}, 65(1):85
  -- 96, 1987.

\bibitem{Hughes1986}
T.~J. Hughes, L.~P. Franca, and M.~Balestra.
\newblock A new finite element formulation for computational fluid dynamics: V.
  {C}ircumventing the {B}abuÃ…Â¡ka-{B}rezzi condition: a stable
  {P}etrov-{G}alerkin formulation of the {S}tokes problem accommodating
  equal-order interpolations.
\newblock {\em Computer Methods in Applied Mechanics and Engineering}, 59(1):85
  -- 99, 1986.

\bibitem{Hughes1989}
T.~J. Hughes, L.~P. Franca, and G.~M. Hulbert.
\newblock A new finite element formulation for computational fluid dynamics:
  {VIII}. {T}he {G}alerkin/least-squares method for advective-diffusive
  equations.
\newblock {\em Computer Methods in Applied Mechanics and Engineering},
  73(2):173 -- 189, 1989.

\bibitem{Hughes1979(1)}
T.~J. Hughes, W.~K. Liu, and A.~Brooks.
\newblock Finite element analysis of incompressible viscous flows by the
  penalty function formulation.
\newblock {\em Journal of Computational Physics}, 30(1):1 -- 60, 1979.

\bibitem{Hughes1978}
T.~J.~R. Hughes.
\newblock A simple scheme for developing ‘upwind’ finite elements.
\newblock {\em International Journal for Numerical Methods in Engineering},
  12(9):1359--1365, 1978.

\bibitem{Johnson1986}
C.~Johnson and J.~Saranen.
\newblock Streamline diffusion methods for the incompressible {E}uler and
  {N}avier-{S}tokes equations.
\newblock {\em Mathematics of Computation}, 47(175):1--18, 1986.

\bibitem{siddiqui2015}
M.~S.~U. Khalid, T.~Rabbani, I.~Akhtar, N.~Durrani, and M.~Salman~Siddiqui.
\newblock {Reduced-Order Modeling of torque on a Vertical-Axis Wind Turbine at
  varying tip speed ratios}.
\newblock {\em Journal of Computational and Nonlinear Dynamics}, 10(4), 2015.

\bibitem{Lovgren2006}
A.~E. L{\o}vgren.
\newblock {\em A reduced basis method for the steady Navier-Stokes problem, in
  ``Reduced basis modelling of hierarchical flow system"}.
\newblock PhD thesis, Norwegian University of Science and Technology, 2006.

\bibitem{Lov2006}
A.~E. Lovgren, Y.~Maday, and E.~M. Ronquist.
\newblock A reduced basis element method for the steady {S}tokes problem.
\newblock {\em ESAIM: Mathematical Modelling and Numerical Analysis},
  40:529--552, 2006.

\bibitem{Manzoni2014}
A.~Manzoni.
\newblock An efficient computational framework for reduced basis approximation
  and a posteriori error estimation of parametrized {N}avier-{S}tokes flows.
\newblock {\em ESAIM: Mathematical Modelling and Numerical Analysis},
  48(4):1199--1226, 2014.

\bibitem{NegriManzoniRozza2015}
F.~Negri, A.~Manzoni, and G.~Rozza.
\newblock Reduced basis approximation of parametrized optimal flow control
  problems for the {S}tokes equations.
\newblock {\em Computers and Mathematics with Applications}, 69(4):319--336,
  2015.

\bibitem{PR014a}
P.~Pacciarini and G.~Rozza.
\newblock Stabilized reduced basis method for parametrized advection-diffusion
  {PDEs}.
\newblock {\em Computer Methods in Applied Mechanics and Engineering},
  274:1--18, 2014.

\bibitem{QV}
A.~Quarteroni and A.~Valli.
\newblock {\em Numerical approximation of partial differential equations},
  volume~23.
\newblock Springer Science \& Business Media, 2008.

\bibitem{Roache1972}
P.~J. Roache.
\newblock {\em Computational Fluid Dynamics}.
\newblock Hermosa Publishers, 1976.

\bibitem{Rovas2003}
D.~Rovas.
\newblock {\em Reduced-basis output bound methods for parametrized partial
  differential equations}.
\newblock PhD thesis, Massachusetts Institute of Technology, 2003.

\bibitem{Rozza2005a}
G.~Rozza.
\newblock {\em Shape design by optimal flow control and reduced basis
  techniques: applications to bypass configurations in haemodynamics}.
\newblock PhD thesis, \'Ecole Polytechnique F\'ed\'erale de Lausanne, N. 3400,
  2005.

\bibitem{Rozza2013}
G.~Rozza, D.~B.~P. Huynh, and A.~Manzoni.
\newblock Reduced basis approximation and a posteriori error estimation for
  {S}tokes flows in parametrized geometries: roles of the inf-sup stability
  constants.
\newblock {\em Numerische Mathematik}, 125(1):115--152, 2013.

\bibitem{Veroy2007}
G.~Rozza and K.~Veroy.
\newblock On the stability of the reduced basis method for {S}tokes equations
  in parametrized domains.
\newblock {\em Computer Methods in Applied Mechanics and Engineering},
  196:1244--1260, 2007.

\bibitem{StabileBallari2018}
G.~Stabile, F.~Ballarin, G.~Zuccarino, and G.~Rozza.
\newblock A reduced order variational multiscale approach for turbulent flows.
\newblock {\em Advances in Computational Mathematics}, 45(5):2349--2368, 2019.

\bibitem{TEZDUYAR1992}
T.~Tezduyar, S.~Mittal, S.~Ray, and R.~Shih.
\newblock Incompressible flow computations with stabilized bilinear and linear
  equal-order-interpolation velocity-pressure elements.
\newblock {\em Computer Methods in Applied Mechanics and Engineering},
  95(2):221 -- 242, 1992.

\bibitem{Veroy2005}
K.~Veroy and A.~T. Patera.
\newblock Certified real-time solution of the parametrized steady
  incompressible {N}avier-{S}tokes equations: rigorous reduced-basis a
  posteriori error bounds.
\newblock {\em International Journal for Numerical Methods in Fluids},
  47(8-9):773--788, 2005.

\end{thebibliography}

\end{document}